\pdfoutput=1
\documentclass[aps,pra,onecolumn,longbibliography,10pt]{article}

\usepackage{subfig}
\usepackage{booktabs}
\usepackage[pdftex]{hyperref}
\usepackage{graphicx}
\usepackage{amsmath}
\usepackage[utf8]{inputenc}
\usepackage{amssymb}
\usepackage{dcolumn}
\usepackage{gensymb} 
\usepackage[normalem]{ulem}

\newcommand{\integer}{\mathop{\rm int}\nolimits}

\usepackage{color}
\title{The Dynamics of Digits: Calculating Pi with Galperin's Billiards}

\author{Xabier M. Aretxabaleta$^{1}$, Marina Gonchenko$^{2}$, Nathan L. Harshman$^{3,}$, Steven Glenn Jackson$^{4}$,\\ Maxim Olshanii$^{5}$ and Grigory E. Astrakharchik$^{1}$\\
{\small $^{1}$Departament de F\'{\i}sica, Universitat Polit\`{e}cnica de Catalunya, 08034 Barcelona, Spain;}\\ {\small xabier.mendez@ehu.eus (X.M.A.); grigori.astrakharchik@upc.edu (G.E.A.)}\\
{\small $^{2}$Departament de Matem\`{a}tiques, Universitat Polit\`{e}cnica de Catalunya, 08028 Barcelona, Spain; mgonchenko@gmail.com}\\
{\small	$^{3}$Department of Physics, American University, 4400 Massachusetts Ave.\ NW, Washington, DC 20016, USA}\\
{\small	$^{4}$Department of Mathematics, University of Massachusetts Boston, Boston, MA 02125, USA; Steven.Jackson@umb.edu}\\
{\small $^{5}$Department of Physics, University of Massachusetts Boston, Boston, MA 02125, USA; Maxim.Olchanyi@umb.edu}}

\date{April 3rd, 2020}

\begin{document}

\maketitle 

\begin{abstract}
	In Galperin billiards, two balls colliding with a hard wall form an analog calculator for the digits of the number $\pi$. This classical, one-dimensional three-body system (counting the hard wall) calculates the digits of $\pi$ in a base determined by the ratio of the masses of the two particles. This base can be any integer, but it can also be an irrational number, or even the base can be $\pi$ itself. This article reviews previous results for Galperin billiards and then pushes these results farther. We provide a complete explicit solution for the balls' positions and velocities as a function of the collision number and time. We demonstrate that Galperin billiard can be mapped onto a two-particle Calogero-type model. We identify a second dynamical invariant for any mass ratio that provides integrability for the system, and for a sequence of specific mass ratios we identify a third dynamical invariant that establishes superintegrability. Integrability allows us to derive some new exact results for trajectories, and we apply these solutions to analyze the systematic errors that occur in calculating the digits of $\pi$ with Galperin billiards, including curious cases with irrational number bases.
\end{abstract}

\section{Introduction}

\begin{figure}[h]
	\centering
	\includegraphics[width=0.4\columnwidth]{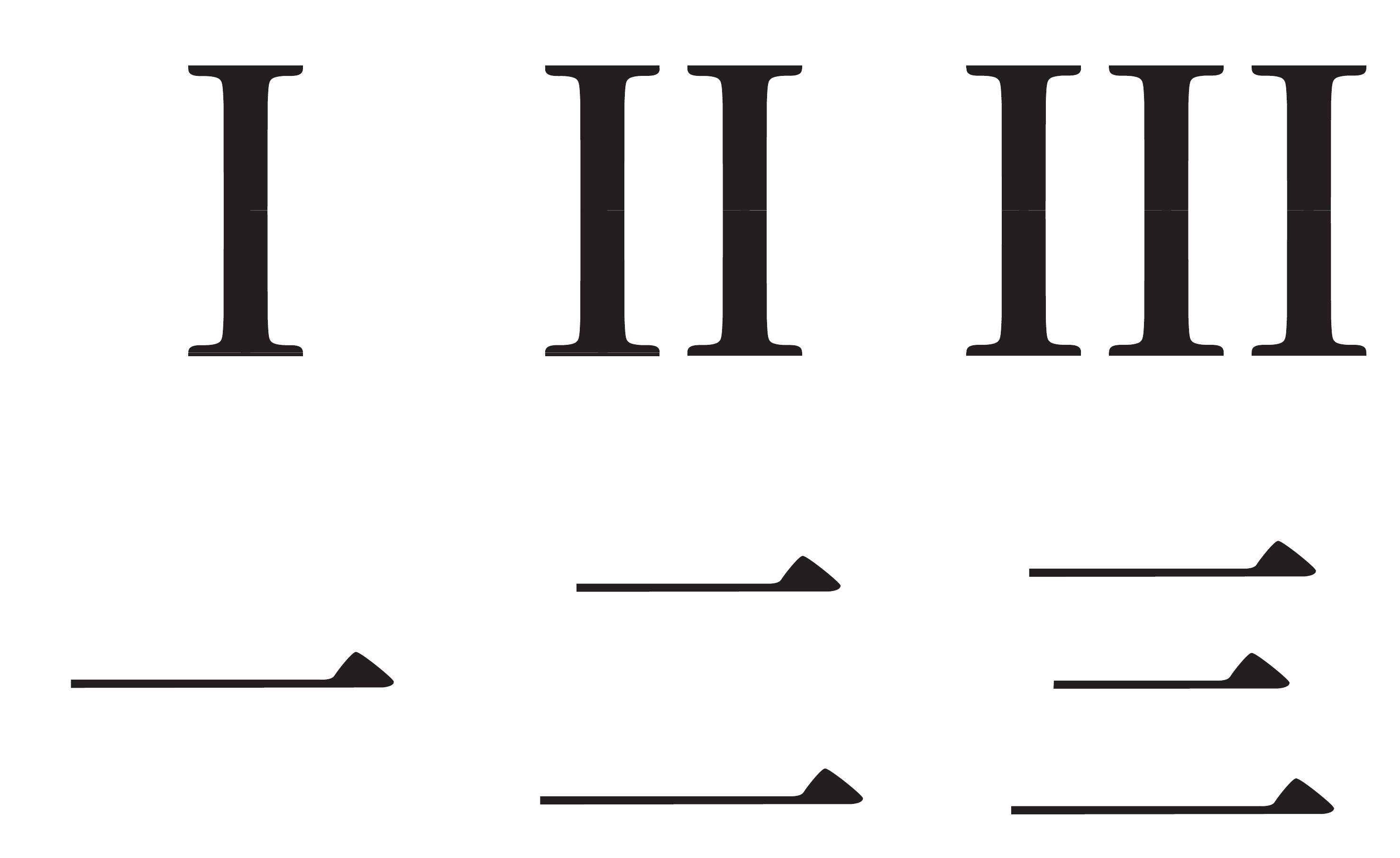}
	\caption{Numbers represented in the unary base.}
	\label{Fig:numerals}
\end{figure}

In the history of science, the invention of numbers stands out as an influential discovery essential for the foundation and development of mathematics and every quantitative science that followed. In many ancient cultures, the symbols for the first digits correspond to a graphical representation of counting. In Babylonian, Roman and Japanese numerals, digit ``1'' contains one counting object, digit ``2'' two objects, digit ``3'' three objects, see Figure~\ref{Fig:numerals}. However, as a technology for representing numbers, the unary base clearly has unfavorable scaling properties. A numeral system, where the position of a digit defines its value with respect to a base that is larger than one, achieves exponential compression at the cost of introducing new symbolic digits. Throughout history different bases have been used, including the modern decimal system and the sexagesimal one introduced in Babylon around the second millennium BC. Its legacy can still be found in modern units of time, with 60 s in one minute and 60 min in one hour.

In a numeral system the digits that represent a rational number eventually repeat, but for irrational numbers like $\pi$ the infinite non-repeating sequence of numbers can be difficult to calculate precisely. Already in the antiquity there was a practical interest in representing numbers like $\sqrt{2}$ and $\pi$ explicitly~\cite{ArndtBook,Goyanes2007,BerggrenBook}. In a Babylonian clay tablet from second millennium BC, the first four digits of $\sqrt{2}$ are explicitly given in sexagesimal system as $1, 24, 51, 10$. In a decimal system the error appears in the eighth digit as can be appreciated by comparing $1 + 24/60 +51/60^2+10/60^3 = 1.4142130$ with the proper value $1.4142135\dots$. The first digit of the number $\pi$, which naturally appears when calculating the ratio between a circumference of a circle and its diameter~\cite{BeckmannBook}, is calculated in the Old Testament \mbox{[1 Kings 7:23]} as $3$. While in many practical situations, it is sufficient to use an approximate value, it~was of fundamental interest to find a method of finding the next digits. Some other ancient estimations come from an Egyptian papyrus which implies $\pi = 256/81 = 3.160\dots$ and a Babylonian clay tablet leading to the value $25/8 = 3.125\dots$. Archimedes calculated the upper bound as $22/7 = 3.1428\dots$.

The fascination with the number $\pi$ makes scientists compete for the largest number of digits calculated. Simon Newcomb (1835-1909) is quoted for having said ``{\it Ten decimal places [of $\pi$] are sufficient to give the circumference of the earth to a fraction of an inch, and thirty decimals would give the circumference of the visible universe to a quantity imperceptible to the most powerful microscope}''\cite{Newcomb1881}. The current world record~\cite{Trueb2016} consists in calculating first 22,459,157,718,361 ($\pi^e$ trillion) digits. Such a task manifestly goes beyond any practical purpose but can be justified by the universal attraction of the number $\pi$ itself.
The distribution of digits is flat in different bases~\cite{Trueb2016} and it was tested that the sequence of $\pi$ digits makes a good random number generator which can be used for practical scientific and engineering computations~\cite{SHU-JU2005}. An alternative popular idea is that, in contrast, special information might be coded in the digits of $\pi$~\cite{ShumikhinBook}, or even God's name as in the plot of ``Pi'' film from 1990. Recently, analogies~between anomalies in the cosmic microwave background and patterns in the digits of $\pi$ were pointed out in ``Pi in the Sky'' article~\cite{Frolop2016}, which appeared on the 1st of April.

While the number $\pi$ elegantly arises in a large variety of trigonometric relations, integrals, series, products, continued fractions, far fewer experimental methods are known of how to obtain its digits by performing measurements according to some procedure. A stochastic method, stemming from Comte de Buffon dates back to the eighteenth century and consists in dropping $N$ needles of length $l$ on parallel lines separated by length $L$ and experimentally determining the number of times $N_{cross}$ that those needles were crossing the lines. The number $\pi$ can be then approximated by $\pi \approx 2l\cdot N/(L N_{cross})$ with the error in its estimation proportional to $1/\sqrt{N}$.
It means that in order to get the first $D$ digits right, one has to perform more than $100^{D}$ trials. This makes it extremely difficult to obtain the precise value in a real-world experiment although an equivalent computer experiment can be easily performed with modern computational power.

\begin{figure}[h]
\centering
\includegraphics[width=0.5\columnwidth]{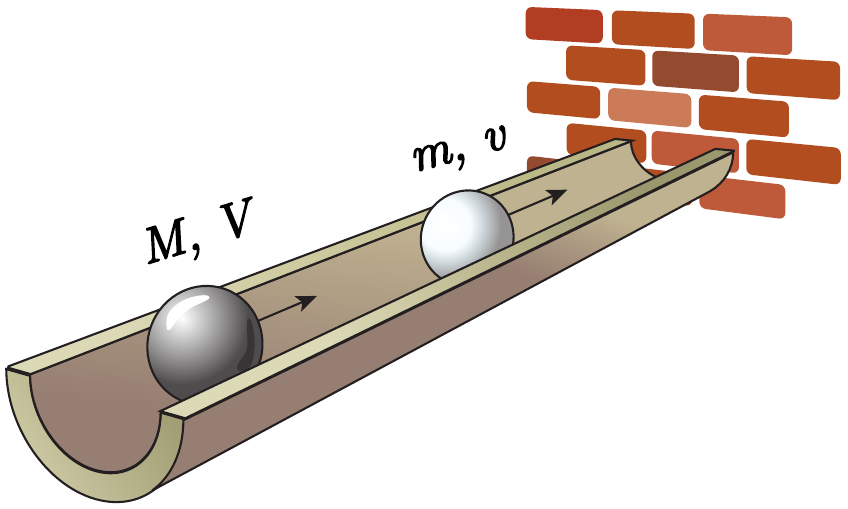}
\caption{Schematic picture of a billiard system, consisting of a heavy ball $M$, light ball $m$ and a wall.}
\label{Fig:balls}
\end{figure}

A completely new perspective has emerged when G.A.~Galperin formulated a deterministic method based on a two-ball billiard~\cite{Galperin2001}. The scheme of the method is summarized in Figure~\ref{Fig:balls}. Two~balls, heavy and light, move along a groove which ends with a wall. The heavy ball collides with the stationary light ball and the number of collisions $\Pi$ is counted for a different mass ratio of the heavy to light ball. It was shown by Galperin that the number of collisions is intimately related to the number $\pi$, providing the first digits of the irrational number. Thus, for equal masses, $M = m$, the~number of collisions is $\Pi = 3$, which corresponds to the first digit of number $\pi$. For masses $M = 100 m$ the number of collisions is $\Pi = 31$, giving the first two digits. The case of $M = 10 000 m$ results in $\Pi = 314$ thus providing three digits and so on. To a certain extent, finding digits of number $\pi$ became conceptually as simple and elegant as enumerating the counting objects like Roman or Japanese digits shown in Figure~\ref{Fig:numerals}.

The history of this elegant method begins with the problem posed by Sinai~\cite{SinaiBook} of the ergodic motion of two balls within two walls. He showed that the configuration space of the system is limited to a triangle and, thus the problem is equivalent to a billiard with the same opening angle. Furthermore, Sinai used ``billiard variables'' such that the absolute value of the rescaled velocity is conserved and the product of the vector of the rescaled velocities with the vector $(\sqrt{M},\sqrt{m})$ is constant. The number of collisions in a ``gas of two molecules'' was given in the book by Galperin and Zemlyakov~\cite{GalperinBook} in 1990, although no relation to the digits of $\pi$ was given at that moment. Similarly, Tabachnikov in his 1995 book~\cite{TabachnikovBook} argued that the number of collisions is always finite and provided the same expression for it. The way to extract digits of $\pi$ from the billiard was explained in Galperin's seminars in the 1990s. In~2001 Galperin published a short article on that procedure in Russian\cite{Galperin2001} and in 2003 in English~\cite{Galperin2003}. This fascinating problem was given as a motivating example in the introduction of another book by Tabachnikov, Ref.~\cite{TabachnikovBook2005}, illustrating trajectory unfolding. Gorelyshev in Ref.~\cite{Gorelyshev2006} applied adiabatic approximation to the two-ball problem and found a conserved quantity, namely the action, close to the point of return. Weidman~\cite{Weidman2013} found two invariants of motion, corresponding to ball--ball and ball--wall collisions and he noted that the terminal collision distinguishes between even and odd digits. Davis in Ref.~\cite{Davis2015} solved the equations of motion as a system of two linear equations for ball--ball and ball--wall collisions, finding the rotation angle from determinantal equation. In addition to the energetic circle~\cite{Galperin2001}, defining the velocities, he expressed the balls positions as a function of the number of collisions.
Related systems studied by similar methods include two balls in one dimension with gravity~\cite{Whelan1990}, dynamics of polygonal billiards~\cite{Gutkin1996} and a ping-pong ball between two cannonballs~\cite{Redner2004}.

In the present work, we review how Galperin billiards with mass ratio $M / m = b^{2N}$ generate the first $N$ digits of the fractional part (i.e., digits beyond the radix point) of number $\pi$ in base-$b$ numeral system. Because Galperin billiards turn the calculation of the digits of an irrational number into a physical process, this motivates deeper inquiry both into the dynamics of this few-body system and into the systematic errors inherent in the analog calculation process. We consider the cases of integer and non-integer base systems, including a compelling case of representing number $\pi$ in a system base of $\pi$.

The main new results of the present paper include the following. We provide a complete explicit solution for the balls' positions, velocities and collision moment as a function of the collision number. We find new invariants of motion and show that for general values of the parameters the system is integrable and for some special values of the parameters it is superintegrable and maximally superintegrable. Finally, we demonstrate that this billiard can be mapped onto a two-particle Calogero-type model.

The article is organized as follows. In Section~\ref{Sec:galperin billiard method}, we review previous results on Galperin billiards. In particular, we explain how billiard coordinates and the unfolding process simplify the analysis of the trajectories in configuration space. We also review prior results from Gorelyshev~\cite{Gorelyshev2006} and Weidman~\cite{Weidman2013} that provide the hints of another dynamical invariant. In Section~\ref{Sec:int}, this invariant is revealed to be a new integral of motion, a kind of pseudo-angular momentum in billiard coordinates. We apply this invariant to generate new analytic results for the equations of motion and to make useful approximations when the mass ratio is large. Using this invariant, we also show that the portrait of the system is close to a circle in velocity-velocity and velocity-inverse position coordinates and to a hyperbola in position-time variables. For certain mass ratios a third integral of motion also exists, making those particular cases of Galperin billiards superintegrable. In Section~\ref{Sec:physical realizations} we discuss different physical systems that can realized Galperin billiards, including finite-size balls, a four-body realization, and we make the connection to Calogero-type models. In Section~\ref{Sec:error} we introduce the concept of systematic error and analyze the possible differences between digits generated by the Galperin billiard method and the usual methods of expressing the number $\pi$ in an arbitrary base. The following Sections~\ref{Sec:integer Bases} and~\ref{Sec:Non-integer Bases} provide examples of how $\pi$ is calculated in systems with integer and non-integer bases, respectively including the intriguing case of expressing $\pi$ in the base $b=\pi$, the generated expression is different from a finite number, $\pi = 1\times \pi^1$, and instead is given by an infinite representation, $\pi = 3 + 1/\pi^2 + 1/\pi^3+\cdots$. The difference between finite and infinite representation is similar to that of $1 = 0.999(9)$ in the decimal system. Furthermore, we note that the finite representation is not unique in the base of the golden number. The conclusion provides a few remarks on how this work can be extended, including to quantum systems where intriguing connections between Galperin billiards and quantum search algorithms have been proposed~\cite{brown}.

\section{Galperin Billiard Method \label{Sec:galperin billiard method}}

In this Section, we summarize known results about the Galperin billiard model and review how it can be used to calculate the digits of $\pi$. The idealized billiard system consists of two `balls' (really, structureless particles) with different masses moving in one dimension and bounded by a hard wall. The initial conditions presuppose that the heavier ball is coming in from infinity and the lighter particle is interposed between the heavy ball and the wall. The ball--ball collisions are perfectly elastic and so are the ball--wall collisions.

Denote the mass and position of the heavier particle by $M$ and $X$ and of the light particle by $m$ and $x$. (As a general rule, we assign capital letters to the variables of the heavy ball and lower-case letters to those of the light ball.)  We choose the coordinate system such that the wall is at the origin and therefore the positions of the heavy and the light balls satisfy $X\leq x \leq 0$ at every moment of time. The heavy ball moves toward the stationary light ball with some initial velocity $V_0>0$ and the light ball begins at rest at the initial position $x_0$. The precise values of $V_0$ and $x_0$ are irrelevant for the total number of collisions, but the position $x_0$ defines the length scale for the system and $x_0/V_0$ sets the time scale.

Somewhat like a force carrier, the light ball shuttles back and forth between the heavy ball and the wall, effectively mediating a repulsive interaction that eventually reverses the approach of the heavy ball. Collisions continue until either:
\begin{enumerate}
\item The last ball--ball collision results in both balls receding from the wall and the heavier ball moving faster. In this case, there are an odd number of collisions.
\item After the last ball--ball collision the heavy ball recedes from the wall with a speed too great for the light ball to catch it again upon one more reflection from the wall. There are an even number of collisions in this case.
\end{enumerate}
Either way, we denote the total number of ball--ball and ball--wall collisions by $\Pi$ and we seek to calculate this number as a function of the mass ratio.
We parameterize the mass ratio as
\begin{eqnarray}
M/m = b^{2N}\;,
\end{eqnarray}
with parameter $b$ (integer or not) referred to as the {\em base} and non-negative integer $N$ referred to as the {\em mantissa}. Using this parameterization for the total number of collisions $\Pi(b,N)$, the integer $N$ determines the number of obtained digits of $\pi$ calculated in the base $b$. 

\subsection{Billiard Coordinates and the Number of Collisions \label{Sec:number of collisions}}

The number of collisions $\Pi$ can be derived most easily by transforming into \mbox{\emph{billiard coordinates}~\cite{SinaiBook,GalperinBook}}. In billiard velocity coordinates the conservation laws implied by the elastic ball--ball and ball--wall collisions take on a simple geometric form~\cite{Galperin2001,Galperin2003,Weidman2013,Davis2015}.

\begin{figure}[h]
\centering
  \includegraphics[width=0.4\linewidth]{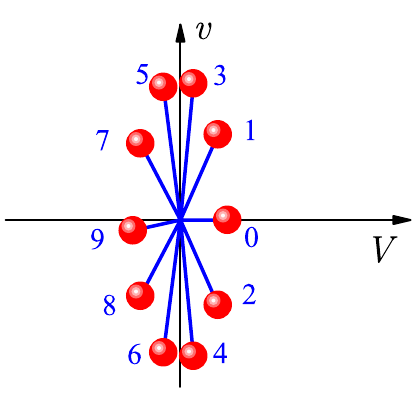}
  \caption{A characteristic example of the dependence of the vector of velocities $\mathbf{v} = (V, v)$ on collision number $n$. The shown data is obtained for $b=3$ and $N=1$ and is an example where $\Pi=9$ and the last collision is ball--ball. The vector of velocities form an ellipse. }\label{Fig:energetic circle:a}
\end{figure}

The conservation of kinetic energy
\begin{eqnarray}
T=\frac{1}{2}M V^2 + \frac{1}{2}m v^2 = \frac{1}{2}M V_0^2 \;
\label{Eq:conservation law:Ekin}
\end{eqnarray}
defines an ellipse in velocity space; see Figure~\ref{Fig:energetic circle:a}. This motivates the change into mass-scaled billiard variables:
\begin{eqnarray}
Y &=& \sqrt{\frac{M}{M+m}} X \equiv \cos(\beta) \ X \nonumber\\
y &=& \sqrt{\frac{m}{M+m}} x \equiv \sin(\beta) \ x.
\label{Eq:billiard variables}
\end{eqnarray}

Note that unlike Sinai or Galperin, we normalize the billiard coordinates by the square root of the total mass so they continue to have units of position. The transformation is described by the angle
\begin{equation}
\tan\beta = \sqrt{m/M} = b^{-N}\label{Eq:beta}.
\end{equation}

The billiard velocities (or configuration speed in Ref.~\cite{GalperinBook}) are defined as time derivative of the position~(\ref{Eq:billiard variables}) and are also scaled with the masses of the balls
\begin{eqnarray}
\nonumber
W&=&\frac{dY}{dt}=\cos\beta \ \frac{dX}{dt}=\cos(\beta) \  V \\
w&=&\frac{dy}{dt}=\sin\beta \ \frac{dx}{dt}= \sin(\beta) \  v.
\label{Eq:billiard velocities}
\end{eqnarray}

The energy conservation law~(\ref{Eq:conservation law:Ekin}) expressed in billiard velocities~(\ref{Eq:billiard velocities}) reads as
\begin{equation}
T = \frac{M+m}{2}\left(W^2+ w^2 \right)  = \frac{M+m}{2}W_0^2.
\label{Eq:conservation law:circle}
\end{equation}

\begin{figure}[h]
\centering
  \includegraphics[width=0.4\linewidth]{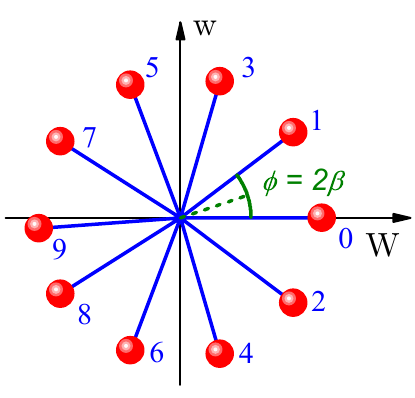}
  \caption{Same scenario as Figure~\ref{Fig:energetic circle:a}, but now the vector of rescaled velocities $\mathbf{w} = (W, w)$ form a circle.}\label{Fig:energetic circle:b}
\end{figure}

The allowed values of $W$ and $w$ forming a circle in billiard velocity space with radius \mbox{$W_0 = \cos\beta \  V_0$}; see Figure~\ref{Fig:energetic circle:b}. Note that Equation~(\ref{Eq:conservation law:circle}) looks like the kinetic energy of a particle in two dimensions with total mass $(M+m)$. Defining the vector of billiard velocities $\mathbf{w}=(W,w)$, the conservation of kinetic energy implies that neither ball--ball nor ball--wall collisions change the magnitude of $|\mathbf{w}|$. Instead, both types of collisions only change the angle $\phi$ between $\mathbf{w}$ and the horizontal axis, as will be explained in more details below.

\unskip

Ball--ball collisions act like reflections on $\mathbf{w}$ in billiard velocity space. To show this, first note that in ball--ball elastic collisions, momentum conservation
\begin{equation}
 p = M V + mv = (M+m)\left(\cos(\beta)\ W + \sin(\beta)\ w \right)  = MV_0,
 \label{Eq:conservation law:momentum}
\end{equation}
combined with energy conservation implies that the relative speed $u$ is also a conserved quantity
\begin{equation}
u = |V - v| = \left|\frac{W}{\cos\beta} - \frac{w}{\sin\beta}\right| = V_0.
\end{equation}

This can be verified by multiplying Equation~(\ref{Eq:conservation law:momentum}) by $V/2$ and subtracting Equation~(\ref{Eq:conservation law:Ekin}). In billiard velocity space, conserved quantities $p$ and $u$ define an orthogonal coordinate system with unit vectors
\begin{eqnarray}
\hat{\mathbf{p}} &=& \begin{pmatrix} \cos\beta \\ \sin\beta \end{pmatrix}\nonumber\\
\hat{\mathbf{u}} &=& \begin{pmatrix} -\sin\beta \\ \cos\beta \end{pmatrix}.
\end{eqnarray}

The coordinate system defined by $\hat{\mathbf{p}}$ and $\hat{\mathbf{u}}$ is rotated by an angle $\beta$
from the $(W,w)$ coordinates. In a ball--ball collision, momentum conservation implies the projection $\mathbf{w}\cdot\hat{\mathbf{p}}$ of the billiard velocity on the momentum axis  is invariant and the projection $\mathbf{w}\cdot\hat{\mathbf{u}}$ on the relative velocity axis  is reversed in sign. In other words, each ball--ball collision acts on $\mathbf{w}$ like a reflection across a line making an angle $\beta$ with the $W$-axis. This can be represented as a matrix acting in billiard velocity space 
\begin{equation}
S_{BB} = \left(\begin{array}{cc} \cos(2 \beta) & \sin(2 \beta) \\ \sin(2 \beta) & -\cos(2 \beta) \end{array}\right),
\end{equation}
or equivalently each ball--ball collision maps $\phi$ to $-\phi + 2\beta$.

Ball--wall collisions are also reflections in billiard velocity space, reversing $v$ (and therefore $w$) while leaving $V$ and $W$ unchanged. Ball--wall collisions are represented as the matrix that preserves the horizontal component of $\mathbf{w}$ and reflects the vertical
\begin{equation}
S_{BW} = \left(\begin{array}{cc} 1 & 0 \\ 0 & -1 \end{array}\right)
\end{equation}

This corresponds to the transformation $\phi$ to $-\phi$. The product $S_{BW}S_{BB}$ of the two reflections is a rotation, and so the new velocities after a combined round of ball--ball and ball--wall collisions are represented in billiard velocity space by a rotation by the angle $2\beta$~\cite{Davis2015}, as depicted in Figure~\ref{Fig:energetic circle:b}.

Using the transformation rules describing ball--ball and ball--wall collisions, it is straightforward to obtain the angle $\phi_{n}$ describing the velocities after $n$ collisions. Typical changes of the vector $\mathbf{w}$ are depicted in Figure~\ref{Fig:energetic circle:b} and can be summarized as follows:
\begin{itemize}
\item $n=0$: Before any collision has happened, the light particle is at rest, $\mathbf{w}_0 = (W_0,0)$, as shown by the horizontal vector with $\phi_0 = 0$
\item $n=1$: The first ball--ball collision reflects the vector $\mathbf{w}_0$ across the line $\phi=\beta$, resulting in $\mathbf{w}_1 = S_{BB} \mathbf{w}_0$ with $\phi_1 = 2\beta$
\item $n=2$: The first ball--wall collision reflects the vector $\mathbf{w}_1$ vertically, resulting in $\mathbf{w}_2 =  S_{BW}\mathbf{w}_1$ with $\phi_2 = -2\beta$
\item $n=3$: The second ball--ball collision reflects the vector $\mathbf{w}_2$ across the line $\phi=\beta$ again, resulting in $\mathbf{w}_3 = S_{BB} \mathbf{w}_2$ with  $\phi_3 = 4\beta$
\item $n=4$: The second ball--wall collision reflects the vector $\mathbf{w}_3$ vertically again, resulting in $\mathbf{w}_4 =  S_{BW}\mathbf{w}_3$ with $\phi_4 = -4\beta$
\end{itemize}

Generalizing, when $n$ is odd, the $n$-th collision is a ball--ball collision after which $\phi_n = (n+1)\beta$. When $n$ is even, the $n$-th collision is a ball--wall collision after which $\phi_n = -n \beta$.

\begin{figure}[h]
\centering
	\includegraphics[width=0.6\columnwidth]{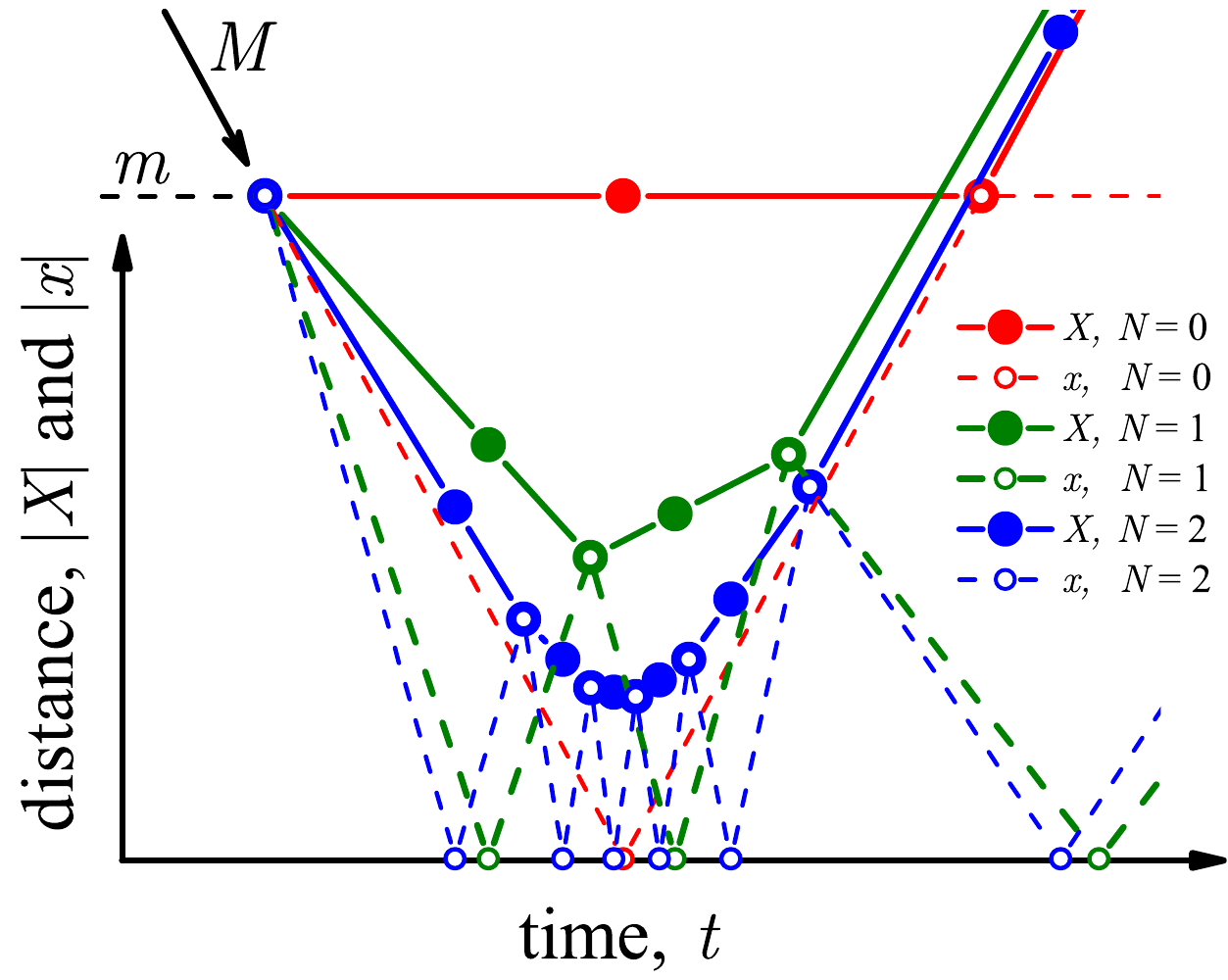}
	\caption{Distance of the heavy and light balls from the wall as a function of time for base $b=2$ and different values of $N$ (in arbitrary units).
		Solid lines and solid symbols, heavy ball $X$;
		dashed lines and open symbols, light ball $x$.
	}
	\label{Fig:positions}
\end{figure}

Working backward from the $\phi_n$, one can integrate the equations of motion to find the times $t_n$ of the collisions and positions $x_n$ and $X_n$ of the balls during the collisions, as for example in \cite{Weidman2013}. See~the Appendix~\ref{Sec:equations of motion} for more details. The last collision defines if the number of collisions is an odd or an even number. Depending on its value, the corresponding digit of $\pi$ is either odd or even. Physically, its~parity depends on whether the last collision was a ball--wall collision with no more ball--ball impacts or if it was a ball--ball collision. In Ref.~\cite{Davis2015} it is shown that an even number of collisions occurs, $\Pi=2k$, when $2k\beta < \pi < (2k+1)\beta$. In Figure~\ref{Fig:positions} we display a typical example of heavy and light ball trajectories, $(t_n, X_n)$ and $(t_n, x_n)$, for $b=2$ and different $N$.

During each ball--ball collision, the velocity $V$ of the heavy ball is changed by a negative amount, eventually stopping and reversing the heavy ball. After the angle $\phi$ has crossed $\pi/2$ position, the~velocity of the heavy mass becomes negative (the ball is moving away from the wall) and its absolute value is increased with each consecutive collision. Collisions continue until $\pi+\beta > \phi_\Pi \geq \pi$ and the lighter ball is at rest or moving away from the wall slower than the heavy ball. After that, the~iterations end; continuing further would result in a decrease of the velocity of the heavy mass, which is physically impossible.

Therefore the ratio of $\pi$ divided by the rotation angle $\beta$ gives the total number of collisions $\Pi$
\begin{equation}
\Pi =  \integer\left[\frac{\pi}{\beta}\right]\;,
\label{Eq:Pi(phi)}
\end{equation}
where $\integer[z]$ means the integer part of $z$.
The number of collisions~(\ref{Eq:Pi(phi)}) can be explicitly evaluated as a function of parameters $b$ and $N$ as
\begin{equation}
\Pi=\integer\left[\frac{\pi}{\arctan(b^{-N})}\right] \;.
\label{Eq:pi:acrctg}
\end{equation}

Moreover, for large base $b$ and large mantissa $N$ the argument of the inverse tangent function is small, $z = b^{-N}\ll 1$, and the inverse tangent function can expanded as $\arctan (z) \approx z$, resulting in Galperin's elegant expression
\begin{equation}\label{Eq:hurbilketa}
\Pi(b,N)\approx \integer\left[\pi b^N \right]
\end{equation}

This equation provides an expression for the number $\pi$ in systems with integer and non-integer bases $b$.

\subsection{Unfolding the Trajectory \label{Sec:Unfolding}}

The trajectory of the balls in configuration space can be given a simple geometrical interpretation which clarifies how the number of collisions is related to the opening angle~\cite{Galperin2003,TabachnikovBook2005} and reveals another invariant of the motion.

\begin{figure}[h]
\centering
\includegraphics[width=0.4\columnwidth]{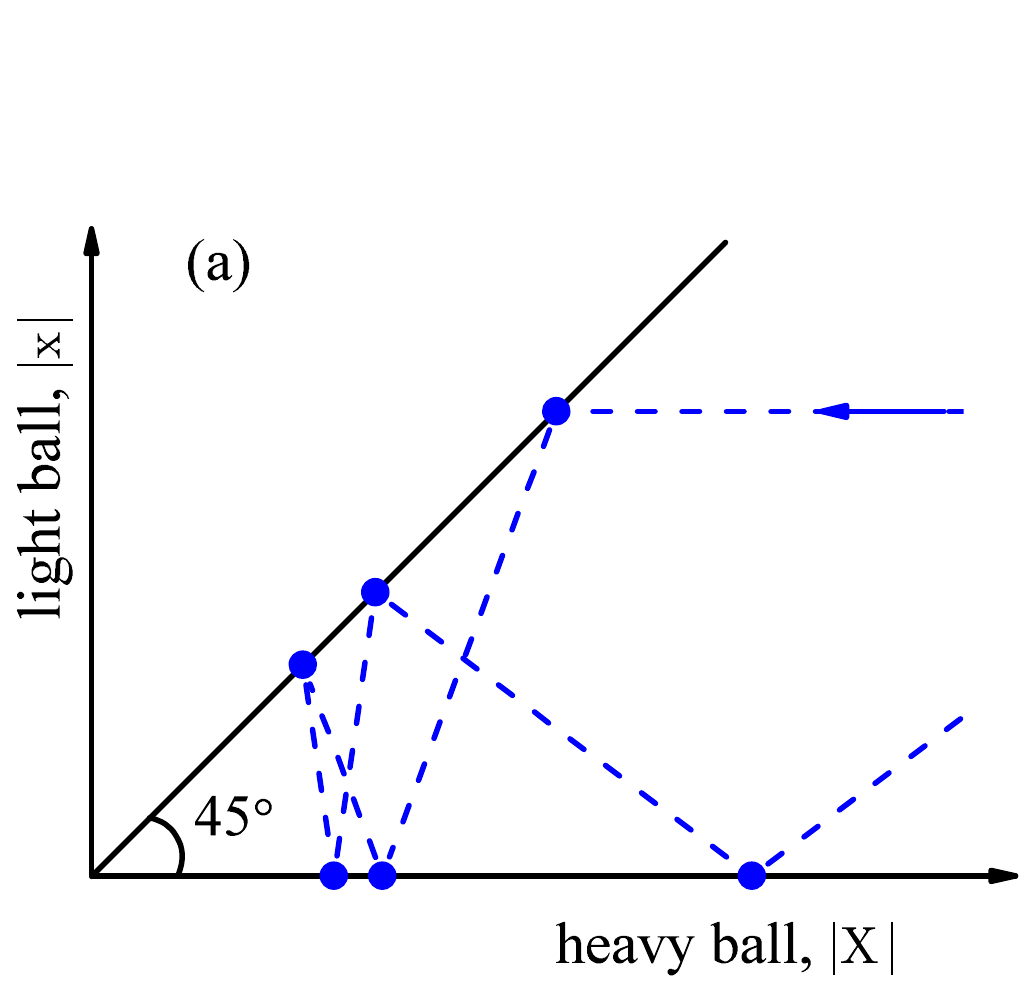}
\includegraphics[width=0.4\columnwidth]{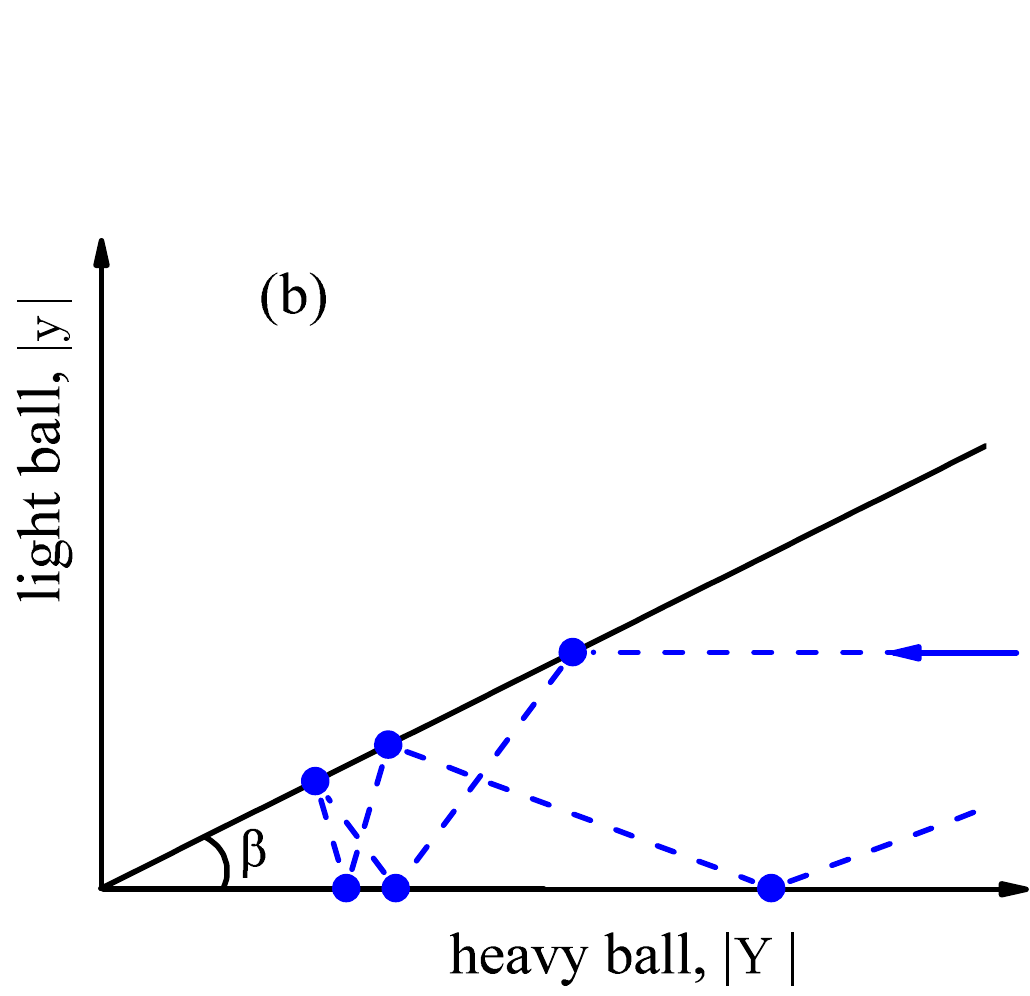}
\includegraphics[width=0.8\columnwidth]{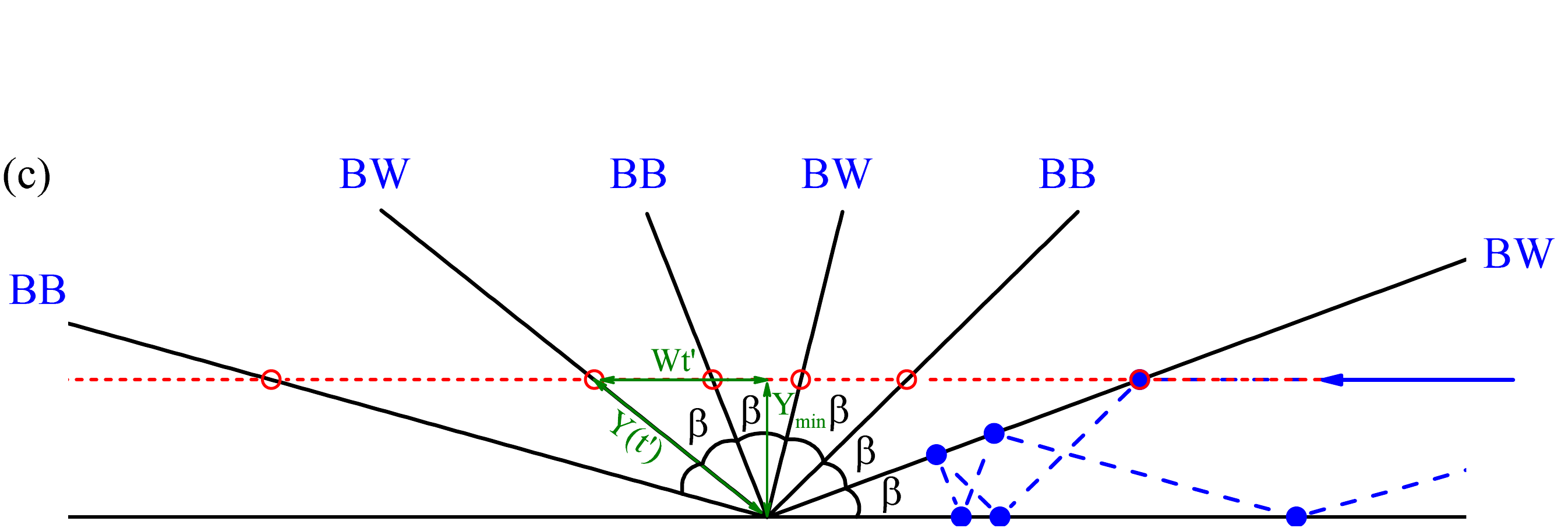}
\caption{{The trajectory in different phase spaces},
(\textbf{a}) original coordinates $0\leq |x| \leq |X|$,
(\textbf{b}) variables of billiard in a wedge, $0\leq |y| \leq |Y|\tan\beta$,
(\textbf{c}) unfolded trajectory.
The parameters are $b=2$ and $N=1$ and correspond to time-dependent data shown in Figure~\ref{Fig:positions}.
}
\label{Fig:unfolding}
\end{figure}

The original particle coordinates are restricted to the region $0\leq |x| \leq |X|$, where boundary $x=0$ corresponds to the light ball hitting the wall and $x=X$ the ball--ball collision, see Figure~\ref{Fig:unfolding}a. The opening angle for this wedge in $(x,X)$ configuration space is $45^\circ$. When the trajectory meets the line $x=0$, corresponding to ball--wall collision, the reflection obeys the law of optics: the tangent component of velocity corresponding to $V$ is preserved and the normal component $v$ is reversed. However, the reflections from the line $x=X$ do not obey the laws of optics, as the incident angle differs from the angle of reflection. 

This is rectified by transforming to billiard coordinates $Y$ and $y$. Now the opening angle in is equal to $\beta$, see Figure~\ref{Fig:unfolding}b and we recover the law of optics for reflections from the line $y=Y \tan\beta$ without disrupting the reflections at $y=0$. This is because the tangent direction to the line $y=Y\tan\beta$ is the unit vector $\hat{\mathbf{p}}$ and the corresponding component of the billiard velocity vector $\mathbf{w}\cdot\hat{\mathbf{p}}$ is proportional to the total momentum and is conserved during ball--ball collisions. Similarly, the normal direction is aligned with the unit vector $\hat{\mathbf{u}}$ and the sign of the normal component is the relative velocity $\mathbf{w}\cdot\hat{\mathbf{u}}$ which is reversed during the collisions.

In this way the original two-body problem is mapped to a problem of a single ball moving in a wedge with opening angle $\beta$ with specular reflections from the mirrors. When the ricocheting trajectory in the wedge is unfolded, it results in a straight line. In Figure~\ref{Fig:unfolding} we show a typical example. It provides a simple geometrical interpretation for the number of collisions as the number of times the opening angle can fit into the maximal angle of $180$ degrees or $\pi$ radians.

The figure of the unfolded trajectory looks like a free particle moving in a straight line, and this is supported by the form of the kinetic energy Equation~(\ref{Eq:conservation law:circle}). This suggests an angular momentum-like quantity
$L= \mathrm{mass} \cdot \mathrm{position} \times \mathrm{velocity}$ should be conserved, at least in magnitude. 
This indeed is the case, and this invariant $L^2= (M+m)^2(Y w - y W)^2$ is the key to the new results presented in Section~\ref{Sec:int}. First, however we present two previous results from Refs.~\cite{Gorelyshev2006,Weidman2013} that anticipated this invariant and provide physical context.

\subsection{Adiabatic Approximation and Action Invariants \label{Sec:adiabatic Approximation}}

From the point of view of Hamiltonian systems, the problem of two balls has two degrees of freedom, namely two positions $X$, $x$ while momenta $P = MV$ and $p=mv$ are conjugate variables. When the heavy ball approaches the point of return near $\phi=\pi/2$, it slows down while the light ball wildly oscillates between the heavy ball and the wall. The heavier the ball, the smaller its minimal distance from the wall, $X_{min}$, and the larger is the maximal velocity of the light ball, $v_{max}$.

This separates the scales into {\em fast} and {\em slow} variables so that during a single oscillation of a light ball, the position of the heavy ball is only slightly changed. It was argued by Kapitza~\cite{Kapitza51UFN,Kapitza51JETP} in his work on driven pendulum (Kapitza pendulum) that by averaging over the fast variables it might be possible to simplify the problem and provide a solution if the separation of scales is large enough. In~our case, the parameter which defines the separation of scales is the mass ratio $M/m = b^{2N}$, so~for any base $b$ by increasing $N$ the needed condition $M/m \gg 1$ is well satisfied. The systems with different scales can be studied in the theory of adiabatic invariants~\cite{ArnoldKozlovNeishtadtBook}.

\begin{figure}[h]
\centering
\includegraphics[width=0.6\columnwidth]{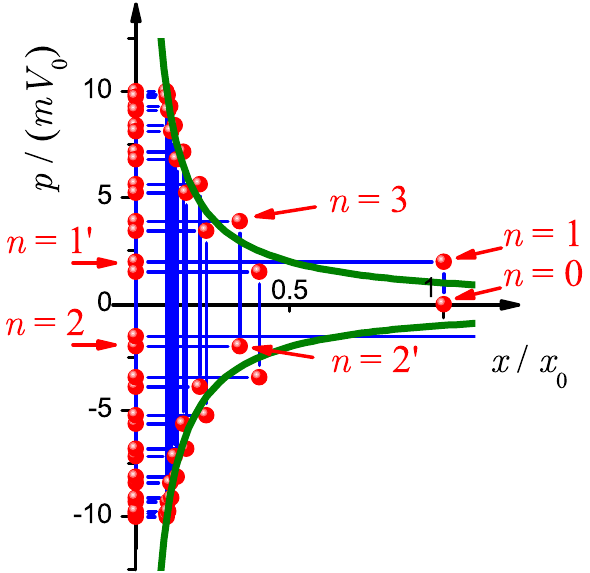}
\caption{(x,p) portrait for $b=10$ and $N=1$;
Red symbols, light particle during ball--wall ($x=0$) and ball--ball ($x\ne 0$) collisions.
Green thick lines, constant action curve defined by Equation~(\ref{Eq:action}).
Blue thin lines, trajectory. 
Index $n=1,2,3,\cdots$ denotes the state after $n$ collisions while primed index $n'$ correspond to an intermediate state in which the velocity of the light ball is not yet reflected.
The area covered by the trajectory between two consecutive collision of the same type (BB or BW) defines the action~(\ref{Eq:action}).
}
\label{Fig:xp}
\end{figure}

It is useful to analyze the $(p,x)$ portrait of the system, corresponding to the fast variables. A~typical example is shown in Figure~\ref{Fig:xp}. After the ball--ball collision, (for example, $n=1$), the light ball moves with a constant momentum $p$ until it hits the wall. This results in a horizontal line with some momentum $p$ and $0<x/x_0 < X_1/x_0$. During the ball--wall ($n=2$) collision the momentum of the light mass is inverted, resulting in a vertical line $x=0$, $p\to -p$. After that the light ball travels with constant momentum until it hits the heavy ball ($n=3$), corresponding to a horizontal line at $-p$ from $0<x/x_0 < X_3/x_0$. At the next ball--ball collision, the velocity of the light particle inverts the sign but its absolute value is slightly changed due to a small but finite momentum transfer from the heavy ball. During a single cycle (or ``period'') consisting of four collisions the light ball draws an almost closed rectangle. The larger is the mass $M$ of the heavy ball and the smaller is its velocity $V$, the more similar is the trajectory during a cycle to a closed rectangle.

The area covered by the light particle during a cycle 
is called {\em action} $I$, defined as
\begin{eqnarray}
I = \frac{  X \left(p - (-p)\right) }{2\pi} = \frac{ X p}{\pi} \,,
\label{Eq:action}
\end{eqnarray}
where $p$ is the momentum of the light particle and $X$ is its maximal distance from the wall (defined by the position of the heavy particle) during a single cycle.
The action~(\ref{Eq:action}) is an adiabatic invariant and is not changed in the vicinity of the return point. Indeed, it might be observed in Figure~\ref{Fig:xp} that while action~(\ref{Eq:action}) is a good adiabatic invariant close to the return point (shown with thick green line), the first few collisions ($n=1; 3; \cdots$) are quite off.

It is shown in Ref.~\cite{Gorelyshev2006} that for times of the order $\varepsilon^2$, action~(\ref{Eq:action}) is conserved with accuracy $\varepsilon$, where $\varepsilon = \sqrt{m/M} = \tan\beta$ is treated as a small parameter. In the same limit the Hamiltonian can be written as
\begin{eqnarray}
H = \frac{P^2}{2M} + \frac{\pi^2 I^2}{2mX^2}\;.
\label{Eq:adiabatic:H}
\end{eqnarray}

It is possible to find two invariants (for BB and BW collisions), which coincide close to the point of return with adiabatic invariant (action) given by Equation~(\ref{Eq:action})~\cite{Weidman2013}.  It can be straightforwardly verified from~(\ref{Eq:equation of motion:BB})--(\ref{Eq:equation of motion:BW}) that the the following quantities remain constant 
\begin{equation}
X_{2k} v_{2k} = \frac{\pi I}{m} = const\;
\label{Eq:invariant:BW}
\end{equation}
for a cycle that starts and ends with a ball--wall collision, and 
\begin{equation}
X_{2k+1} (V_{2k+1}-v_{2k+1}) = \frac{\pi I}{m} = const\;.
\label{Eq:invariant:BB}
\end{equation}
during a cycle between ball--ball collisions. Importantly, these invariants are {\em always} conserved on the corresponding BW and BB cycles and not only close to the point of return.

The ball--wall invariant~(\ref{Eq:invariant:BW}) reduces to an elegant expression, $X_{2k} v_{2k} = - x_0 V_0$, and thus $I$, entering into Equations~(\ref{Eq:invariant:BW}) and (\ref{Eq:invariant:BB}), can be expressed in terms of the initial conditions as
\begin{equation}
I = \frac{|x_0| V_0 m}{\pi},
\label{Eq:action invariant}
\end{equation}
agreeing with (\ref{Eq:action}) above. Furthermore, from the BW invariant we can get the expression for the closest position $X_{min}$ of the heavy ball to the wall  which is reached at the collision when the heavy ball inverts its velocity and the light ball achieves its maximal velocity: 
\begin{equation}
X_{min}
= \frac{V_0}{v_{max}}x_0
= \sqrt{\frac{m}{M}}x_0
= \frac{x_0}{b^N}\;.
\label{Eq:X:min}
\end{equation}
For the ratio of the velocities we have used the energy conservation law~(\ref{Eq:conservation law:Ekin}), since $v_{max}$ is reached when $V=0$, and $v=0$ is reached for $V=V_0$, c.f~Fig.~\ref{Fig:energetic circle:a}.  

Defining $\alpha$ as the phase angle conjugate to $I$, the time derivative of the phase is obtained from Hamiltonian~(\ref{Eq:adiabatic:H}) as $\dot\alpha = dH / dI$. The integration over the time gives the final phase after all collisions have happened as $\alpha^{\rm final} = \pi^2 \sqrt{M}/\sqrt{m}$~\cite{Gorelyshev2006}. During each cycle there are two collisions (BB and BW) and the phase changes by $2\pi$, so the total number of collisions $\Pi$ can be infer as $\alpha^{\rm final} = \Pi \pi$, resulting~in $\Pi = \pi / \beta + O(\beta)$ where $\beta \approx \sqrt{m/M}$. This formula should be compared with Equation~(\ref{Eq:Pi(phi)}) and, indeed, correctly relates total number of collisions $\Pi$ with $\pi$. At the same time it is not \emph{a priori} obvious that the adiabatic approximation should be precise far from the return point, that is for times $t\gg \varepsilon^2$, especially at the time of the final collision. The BB and BW invariants~(\ref{Eq:invariant:BW})--(\ref{Eq:invariant:BB}) coincide with adiabatic invariant (action) given by Equation~(\ref{Eq:action}) close to the point of return, and, in particular, this clarifies why the adiabatic approximation leads to the correct number of collisions even if the region of applicability of the approximation is violated.

Finally we note that the time dependence of the phase $\alpha(t)$ is related to the time dependence of the collision number $n(t)$ according to $\alpha(t_n) = \pi n(t_n)$. In the continuous limit of many collisions, the phase increases as an inverse tangent function, as shown in Equation~(\ref{Eq:t:approx}) of Section~\ref{Sec:exact Solution} below.

\section{Integrability and Its Consequences \label{Sec:int}}

We show now that the Galperin billiard system is \emph{Liouville integrable}, i.e.,\ it possesses as many exact constants (first integrals) of motion  in involution as it has degrees of freedom. Since this system has only two degrees of freedom, this requires the existence of only one constant of motion in addition to the total energy. Generally, a notion of \emph{involution} between two observables---vanishing of the Poisson bracket between them---requires a Hamiltonian reformulation of the laws of dynamics of the system. However, for the two-dimensional systems, the only zero Poisson bracket required is the one between the additional conserved quantity and the Hamiltonian; the latter is simply automatic.

To identify this additional conserved quantity, let us consider the system in the billiard coordinates~(\ref{Eq:billiard variables}). It is represented by a two-dimensional particle with mass $M+m$ moving a wedge of an opening $\tan\beta = \sqrt{m/M}$ as shown in Figure~\ref{Fig:unfolding}. In between the collisions, the angular momentum, $L=(M+m)(Y w - y W)$, is conserved. Upon a ball--wall or a ball--ball collision, the angular momentum changes sign. However, its square,
\begin{equation}
L^2 = (M+m)^2(Y w - y W)^2 = mM (X v - x V)^2
\label{LSq}
\end{equation}
remains invariant throughout the evolution. 

This can be checked by explicit calculations using the results of Appendix~\ref{Sec:equations of motion}. At the instances of a ball--wall collision, where $x=0$, the invariant~(\ref{LSq}) is proportional to the square of the invariant (\ref{Eq:invariant:BW}) identified by Weidman~\cite{Weidman2013}:
\begin{equation}
L^2 \Big|_{2k} = m M (X_{2k} v_{2k})^2 = \pi^2 \frac{M}{m} I^2\,\,.
\end{equation}

Likewise, on a ball--ball collision, the angular momentum square assumes a value proportional to the corresponding invariant~(\ref{Eq:invariant:BB}):
\begin{equation}
L^2 \Big|_{2k+1} = m M \left(X_{2k+1} (v_{2k+1}-V_{2k+1})\right)^2 = \pi^2 \frac{M}{m} I^2\,\,,
\end{equation}
with the same coefficient of proportionality. Evaluating this constant on the initial conditions, the~invariant takes the value 
\begin{equation}
L^2 = m M x_0^2 V_0^2.
\end{equation}

To express the Hamiltonian in terms of this invariant, we start from the unfolded billiard coordinates. The Hamiltonian has the form of a free particle in two dimensions  with mass $M+m$ (c.f.~\ref{Eq:conservation law:circle}) 
\begin{equation}
H = \frac{1}{2(M+m)} \left(P_Y^2 + P_y^2 \right).
\end{equation}
where the momenta conjugate to the unfolded billiard coordinates $(Y,y)$ are $P_Y=(M+m)W$ and $P_y=(M+m)w$. Changing to polar coordinates
\begin{eqnarray}
r &=& \sqrt{Y^2 + y^2}\nonumber\\
P_r &=& (M+m)\dot{r} = (M+m)\frac{YW+yw}{r}\nonumber\\
\phi &=& \arctan(y/Y)\nonumber\\
L &=& (M+m)r \dot{\phi} = (M+m)\frac{Yw-yW}{r}
\end{eqnarray}
the Hamiltonian take the form
\begin{equation}\label{eq:Calogero}
H = \frac{1}{2(M+m)} \left(P_r^2 + \frac{L^2}{r^2} \right).
\end{equation}
In this form, the Hamiltonian~(\ref{eq:Calogero}) looks like a two-body Calogero-Sutherland model as will be explained in more details in Section~\ref{Sec:CSM}.

\subsection{Exact Solution \label{Sec:exact Solution}}

In the Appendix~\ref{Sec:equations of motion}, the sequence of positions and times of collisions was found by step-wise integration of the kinematics. However, the integrability of the Galperin model and the geometry of the unfolded trajectory allows for explicit solution where the velocities and positions can be explicitly expressed as a function of the collision number.

\begin{figure}[h]
\centering
\includegraphics[width=0.8\columnwidth]{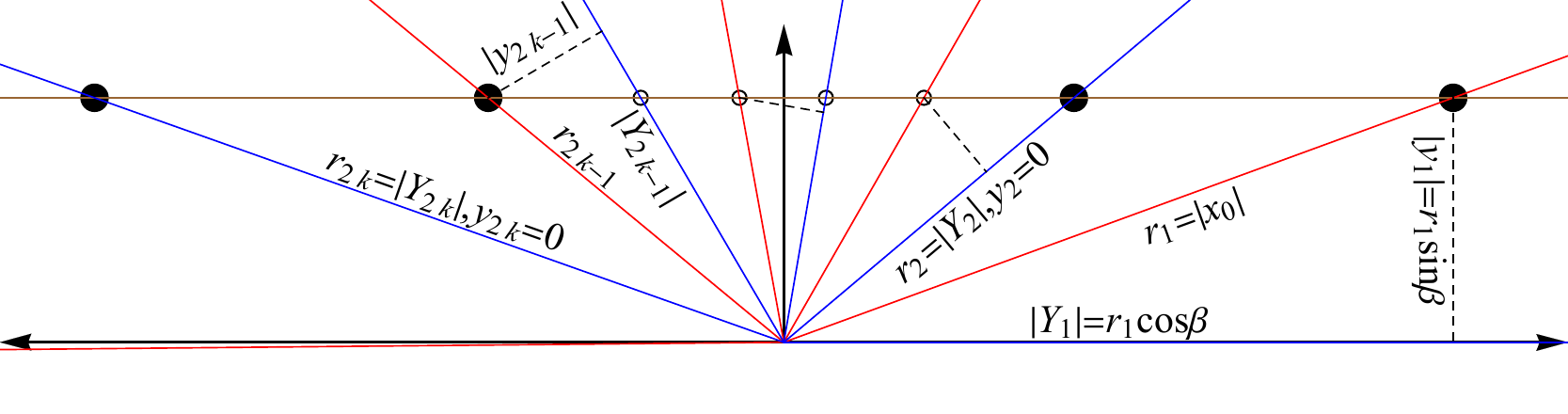}
\caption{Unfolded trajectory depicting relations among trajectory and billiard coordinates. The brown line is the trajectory specified by the initial conditions with $L^2=m M x_0^2 V_0^2$ and speed $W_0=V_0\cos\beta$.  For this mass ratio $M/m = 7.5$, there are eight collisions, depicted as four disks for points with analysis and four circles for other points. Odd-numbered ball--ball collisions occur where the trajectory intersects red lines and even-numbered ball--wall collisions occur at blue lines. The angle for the $n$-th collision point is $n \beta$ and the distance to the origin is $r_n$. The projection of the ball--ball collisions onto the blue ball--wall lines are depicted by dashed lines.}
\label{Fig:geo}
\end{figure}

As shown in Figure~\ref{Fig:geo}, the trajectory in unfolded coordinates is a horizontal line that traverses at constant speed $W_0 = V_0 \cos\beta$. The $n$-th collision occurs at the intersection of the trajectory and a line making an angle $n\beta$, where odd $n$ are ball--ball collisions and even $n$ are ball--wall. These intersections occur at a distance $r_n$ from the origin of the unfolded coordinates, with the first ball--ball collision occurring at $r_1=|x_0|$. Using the law of sines, the general formula for the `collision radius' $r_n$ is
\begin{equation}
    r_n = |x_0| \frac{\sin\beta}{\sin(n\beta)}.
\end{equation}

For ball--ball collisions $n=2k+1$ with $k=0,1,2,\ldots$, we project the collision with radius $r_{2k+1}$ and angle $(2k+1)\beta$ onto the previous ball--wall axis at $2k\beta$ and use (\ref{Eq:billiard variables}) to find
\begin{equation}
x_{2k+1} = X_{2k+1} = -r_{2k+1}= x_0 \frac{\sin\beta}{\sin\left((2k+1)\beta\right)},\ k=0,1,2,\ldots.
\end{equation}

For ball--wall collisions $n=2k$ with $k=1,2,\ldots$, the little ball is at the wall
\begin{equation}
    x_{2k} =0,\ k=1,2,\ldots,
\end{equation}
and the large ball is at
\begin{equation}
    X_{2k} = -\frac{r_{2k}}{\cos\beta} = x_0 \frac{\tan\beta}{\sin(2k\beta)},\ k=1,2,\ldots.
\end{equation}

Following similar geometrical logic, the time interval $t_n$ from the first collision to the $n$-th collision can be found by considering the length of trajectory $d_n$ between those collisions
\begin{equation}
    d_n = |x_0| \frac{\sin(n\beta - \beta)}{\sin(n\beta)}
\end{equation}
and then dividing by the trajectory velocity in billiard coordinates $W_0$ to find
\begin{equation}
    t_n = \frac{|x_0|}{V_0} \frac{\sin(n\beta - \beta)}{\sin(n\beta)\cos\beta} = t_0 \left(1 - \cot(n\beta)\tan\beta\right),\label{eq:totaltime}
\end{equation}
where $t_0 = |x_0|/V_0$ is the characteristic time of the system.
From this the time interval $\tau_n$ between the $(n-1)$-th collision and the $n$-th collision is calculated to be
\begin{equation}
\tau_n = t_n - t_{n-1} = t_0 \frac{\sin\beta\tan\beta}{\sin(n\beta - \beta)}{\sin(n\beta)}.
\label{eq:timeinterval}
\end{equation}

In the continuous limit of many collisions, $b^N\gg 1$, the following simple expression for the inverse of time $\tau_n$ between consecutive collisions,
\begin{eqnarray}
\frac{t_0}{\tau_n}
 \approx b^{2N}\sin^2(n/b^N)
 \approx b^{2N}\cos^2(n'/b^N)
\label{Eq:tau(n):approx}
\end{eqnarray}
where $n' = n - \Pi/2$. In the same limit, we find the following approximate relations between time and collision number
\begin{eqnarray}
\frac{t_{n'}}{t_0} \approx 1 - \frac{1}{b^N}\cot(n/b^N) \approx 1 + \frac{1}{b^N}\tan(n'/b^N).
\label{Eq:t:approx}
\end{eqnarray}

The total time interval from first collision to final collision is 
\begin{equation}\label{eq:totalt}
    t_\Pi = t_0 \left(1 - \cot(\Pi\beta)\tan\beta\right) = t_0 \left(1 - \cot\left(\integer\left[\frac{\pi}{\beta}\right]\beta\right)\tan\beta\right).
\end{equation}

This curious expression for $t_\Pi$ is bounded by below by $2 t_0 = 2|x_0|/V_0 $, i.e.,  the time that the large ball would have taken to transit from $x_0$ to the wall and back to $x_0$ if the small ball \emph{had not been there at all}. This lower bound is saturated in the limit from above when $\beta= \pi/q$ and $q$ is an integer; the~limit taken from below diverges to $t_\Pi\to \infty$. This divergence occurs in the last time step $\tau_\Pi$, as one sees from the expression for $t_{\Pi-1}$ which is bounded from above by $2 t_0 = 2|x_0|/V_0$. Further note that when $\beta= \pi/q$ and $q$ is an integer, the velocity of the large ball is exactly reversed. Both of these effects for $\pi/q$ provide a clue that the system is superintegrable for certain mass ratios, as we show below.

For completeness, we note that the velocities immediately after each collision can be found from the results of Section~\ref{Sec:galperin billiard method} and are tidily expressed as
\begin{eqnarray}
\phi_n &= &(-1)^{n+1}2\beta \integer\left[\frac{n+1}{2}\right], \\
\phi_{2k} &= & - 2k\beta, \;\;\; k=0, 1, \ldots,\\
V_n &= &V_0 \cos \phi_n,\\
v_n &= &V_0 \cot\beta \sin \phi_n.
\end{eqnarray}

\subsection{Position as a Function of Time: Hyperbolic Shape \label{Sec:xvs.t}}

Here we will demonstrate that a hyperbolic curve describes the positions of the light ball at BB collisions and of the heavy ball both at BB and BW collisions.

In the description of a billiard in a wedge, the trajectory is bounded to the phase space $0\leq |y| \leq |Y| \tan\beta$, as shown in Figure~\ref{Fig:unfolding}. The collisions happen when either $y=0$, i.e., when the light particle hits the wall (BW collision) or when $y=Y \tan\beta$ and the light particle hits the heavy one (BB~collision). The unfolded trajectory is formed by reflecting the wedge, so that its angle $\beta$ is preserved. The collisions in the unfolded trajectory occur when the straight line intersects one of the mirrors, corresponding to an angle $n\beta$ with $n$ being the number of the collision. For any intersection its distance from the origin is the same in unfolded picture and that of the billiard in a wedge. In particular, for~a ball--ball collision, this distance is equal to $\sqrt{Y^2(t)+y^2(t)}$. Instead, in the moment of a ball--wall collision, the light ball has coordinate $y(t)=0$ and this distance is equal to the position of the heavy ball $Y(t)$. In the unfolded coordinates depicted in Figure~\ref{Fig:unfolding}), the minimal possible distance $Y_{min}$  of the heavy ball from the wall is located on the vertical line directly above the origin and corresponds to the point of return in the limit $M \gg m$. The projection to the horizontal axis is given by $Wt'$, where $t'$ is the time counted from the point of return and $W$ is the constant velocity, equal to the initial velocity of the heavy ball $W = W_{max}$. Catheti $Y_{min}$, $W_{max}t'$ and hypotenuse $Y(t)$ forming a right-angled triangle are related as $Y^2(t) = Y^2_{min} + (W_{max}t')^2$. The same expression written in terms of the original coordinate $X(t')$ and velocity $V$ leads to the hyperbolic relation
\begin{eqnarray}
\left(\frac{X(t')}{X_{min}}\right)^2 - \left(\frac{t'}{X_{min}/V_{max}}\right)^2 = 1\;,
\label{Eq:hyperbola:BW}
\end{eqnarray}
exactly satisfied for any ball--wall collision. Here $X_{min}$ is given by Equation~(\ref{Eq:X:min}).

Instead, for a ball--ball collision both coordinates of the heavy and light particles are equal, $X=x$, and lie on a hyperbola of a slightly smaller semi-axis
\begin{eqnarray}
\left(\frac{X(t')}{\sqrt{\frac{M}{M+m}}X_{min}}\right)^2 - \left(\frac{t'}{X_{min}/V_{max}}\right)^2 = 1\;,
\label{Eq:hyperbola:BB}
\end{eqnarray}

In the limit of large mass, Equations~(\ref{Eq:hyperbola:BW}) and (\ref{Eq:hyperbola:BB}) coincide.

\begin{figure}[h]
\centering
\includegraphics[width=0.6\columnwidth]{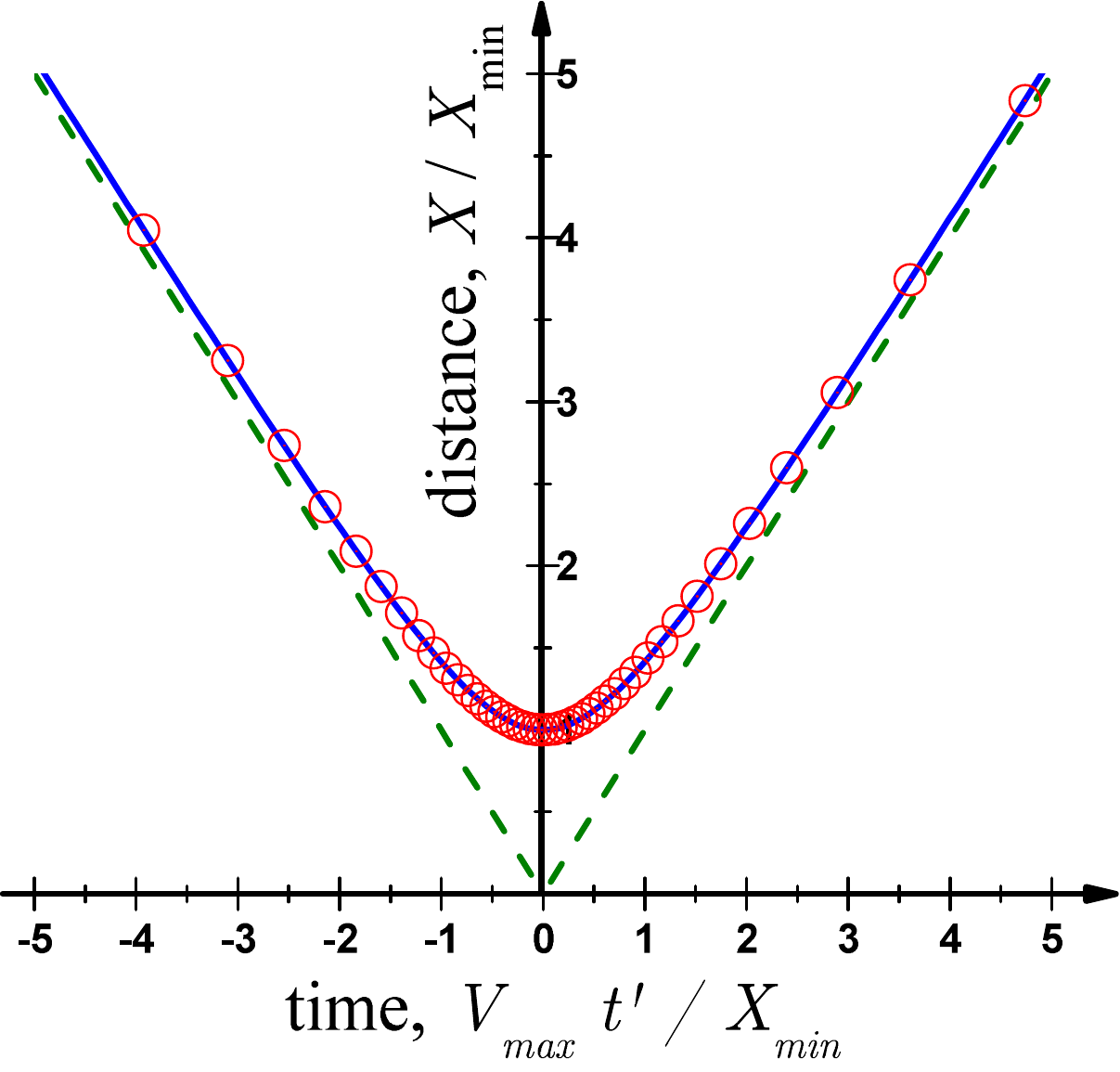}
\caption{Distance of the heavy ball close to the return point, for $b=2$ and $N=3$.
Circles, heavy~ball $X$;
dashed line, limit of an infinitely massive ball, Equation~(\ref{Eq:return:heavyball});
solid line, hyperbola defined by Equation~(\ref{Eq:hyperbola:BB}).
}
\label{Fig:hyperbola}
\end{figure}

We compare predictions of Equation~(\ref{Eq:hyperbola:BB}) with the exact results in Figure~\ref{Fig:hyperbola}. The minimal possible distance $X_{min}$ is actually reached only if there is a crossing of the unfolded trajectory at the vertical line above the origin (see Figure~\ref{Fig:unfolding}), otherwise the actual minimal distance is larger. In the limit of an infinitely heavy ball ($N\to\infty$ and $M\to\infty$), the trajectory again becomes non-analytic with a kink in the point where the heavy ball hits the wall,
\begin{eqnarray}
X(t') = V_{max}|t'|
\label{Eq:return:heavyball}
\end{eqnarray}
shown in Figure~\ref{Fig:hyperbola} with two straight lines.

\subsection{Circle in \texorpdfstring{$(V,1/X)$}{(V,1/X)} Variables \label{Sec:circle}}

Within the adiabatic approximation, introduced in Section~\ref{Sec:adiabatic Approximation}, from the Hamiltonian~(\ref{Eq:adiabatic:H}) it can be seen that the $(P, 1/X)$ portrait has a semicircular shape with the coefficient of proportionality linear in the action $I$~\cite{Gorelyshev2006}. In Section~\ref{Sec:adiabatic Approximation} it was verified that some of the predictions of the adiabatic theory actually remain exact even far from the point of return, effectively expanding the limits of its applicability. In particular, the action $I$ is conserved for any ball--wall collision throughout the whole process, and its value can be expressed in terms of the initial position of the light ball $x_0$ and the initial velocity of the heavy ball $V_0$ according to Equation~(\ref{Eq:action invariant}). This suggests that the portraits in ($P$, $1/X$) and ($V$, $1/X$) coordinates are close to ellipses.

A straightforward way to see it is to use the Hamiltonian~(\ref{Eq:adiabatic:H}) obtained within the adiabatic~approximation,
\begin{eqnarray}
\frac{MV^2}{2} + \frac{\pi^2 I^2}{2mX^2}
= \frac{MV_0^2}{2}
\label{Eq:H:(V,X)}
\end{eqnarray}
where we extended its validity to any BW collision, in particular to collisions happening far from the point of return, $X\to\infty$. Equation~(\ref{Eq:H:(V,X)}) can be recast in the form of an ellipse for $(V,1/X)$ coordinates as
\begin{eqnarray}
\frac{1}{b^{2N}}
\left(\frac{x_0}{X}\right)^2
+
\left(\frac{V}{V_0}\right)^2
=
1\;.
\label{Eq:ellipse:(V,X)}
\end{eqnarray}

\begin{figure}[h]
\centering
\includegraphics[width=0.5\columnwidth]{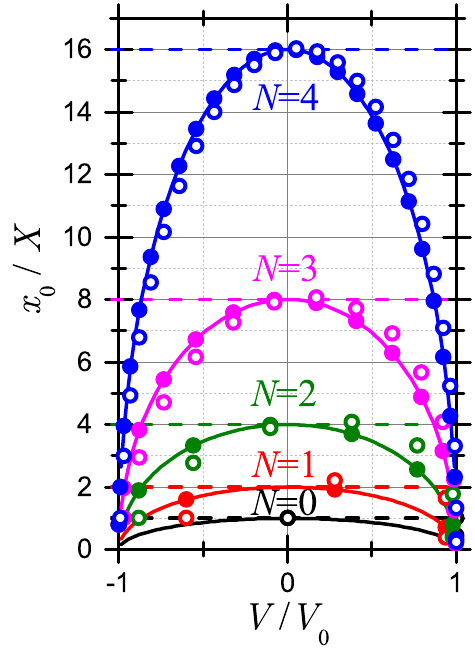}
\caption{
$(V,1/X)$ portrait for $b=2$ and $N=0,1,2,3,4$ (from bottom to top).
Solid symbols, ball--wall collision;
open symbols, ball--ball collision.
Solid lines, ellipse~(\ref{Eq:ellipse:(V,X)}).
Dashed horizontal lines, inverse of the minimal distance~(\ref{Eq:X:min}).
}
\label{Fig:ellipse}
\end{figure}

Figure~\ref{Fig:ellipse} shows an example of the trajectory in $(V,1/X)$ coordinates. The first collision happens at $V/V_0 = 1$  and $x_0/X = 1$ corresponding to the initial velocity and the initial (large) distance from the wall. As the collisions go on, the heavy ball comes closer to the wall until it inverts its velocity at the point $V=0$, which corresponds to the point of return. At this moment the heavy ball is located at the closest distance to the wall. It might be appreciated that Equation~(\ref{Eq:X:min}) describing this distance is quite precise from the practical point of view. The case illustrated in Figure~\ref{Fig:ellipse} corresponds to the binary base, $b=2$, and for mantissa length $N=1;2;3;4$ the heavy ball is expected to come closer to the wall by a factor of $2;4;8;16$ compared to the initial position of the light ball. Once the point of return is passed, the heavy ball has a negative velocity which increases in absolute value up to $V/V_0 \approx -1$, while the ball moves far away from the ball $x_0/X\to 0$.

Overall, the shapes obtained are quite similar to ellipses predicted by Equation~(\ref{Eq:ellipse:(V,X)}). The~``discretization'' becomes smaller as $N$ is increased. The pairs with same velocity $V$ but different values of $X$ correspond to light ball--wall collisions (velocity of the heavy ball is not changed) and ball--ball collisions, shown in Figure~\ref{Fig:ellipse} with closed and open symbols. The points which correspond to the ball--wall collisions lie exactly on the top of the ellipse due to presence of ball--wall invariant~(\ref{Eq:invariant:BW}). Instead, for ball--ball collisions there is some shift, with a different sign for the heavy ball moving towards the wall or away from it. This effect originates from an additional contribution containing the velocity of the heavy ball in the ball--ball invariant~(\ref{Eq:invariant:BB}). At the point of return this correction vanishes and the adiabatic theory becomes fully applicable.

\subsection{Superintegrability and Maximal Superintegrability \label{Sec:super}}

When a system with $d$ degrees of freedom has more than $d$ constants of the motion, it is called \emph{superintegrable}. The maximum allowed number of functionally-independent conserved quantities is $2d-1$, one less than the dimension of the phase space. Such a system is called \emph{maximally superintegrable}. For certain mass ratios, the Galperin model has a third functionally-independent integral of motion. Since the Galperin model has two degrees of freedom, $d+1 = 2d -1$ and it is therefore superintegrable and also maximally superintegrable.

For a system with bounded orbits, the superintegrability manifests itself as a reduction of the dimensionality of the phase space available from a given initial condition. Maximal superintegrability results in closed one-dimensional orbits. For unbounded orbits like the Galperin model, the manifestation of the maximal superintegrability is more subtle, but still---as shown below---tangible.

For certain mass ratios
\begin{equation}
\frac{m}{M} = \tan^2(\pi/q) 
\label{kaleidoscopes}
\end{equation}
the wedge in the billiard coordinates $(y,Y)$ depicted in Figure~\ref{Fig:unfolding}b acquires an opening of $\beta = \pi/q$ with $q \ge 3$ being an integer. For these rational angles, a third functionally-independent constant of the motion appears. In this case, sequences of reflections about the cavity walls form a finite group with order $2q$ known as the reflection group $I_{2}(q)$ or the dihedral group $D_q$. The generators for the group are the reflections in billiard velocity space $S_{BB}$ and $S_{BW}$ defined in Section~\ref{Sec:number of collisions} supplemented with the group-defining relation $(S_{BW}S_{BB})^q =1$.

As it has been shown in Ref.~\cite{olshanii2015_105005}, this discrete reflection group symmetry implies that a new constant of motion can be constructed: it is represented by the first nontrivial invariant (or Chevalley) polynomial of the group \cite{chevalley1955_778,mehta1988_1083}, evaluated on the momentum vector. The constant of motion $J$ produced by this construction in our case is as follows:
\begin{equation}
J
= \frac{1}{2 \cos^q\beta}\left((W+iw)^{q}+(W-iw)^{q}\right)
= \frac{1}{2}
\left(
   (V + i \tan(\pi/q) v)^{q} + (V - i \tan(\pi/q) v)^{q}
\right)
\label{J}
\,\,.
\end{equation}

Some notable examples include
\begin{align*}
&
q=3
&
\frac{m}{M} = 3
&&
J = V^3-9 V v^2
&
\\
&
q=4
&
\frac{m}{M} = 1
&&
J = V^4 - 6 V^2 v^2 +v^4
&
\\
&
q=5
&
\frac{m}{M} = 5 - 2 \sqrt{5}
&&
J = V^5 - 10(5-2\sqrt{5}) V^3 v^2 + 25(9-4\sqrt{5}) V v^4
&
\\
&
q=6
&
\frac{m}{M} = \frac{1}{3}
&&
J = V^6 - 5 V^4 v^2 + \frac{5}{3} V^2 v^4 - \frac{1}{27} v^6
&
\,\,.
\end{align*}

Observe that in these examples and in general, even (odd) $q$, produces a constant of motion $J$, which is an even (odd) function with respect to the $V\to -V,\,v\to -v$ inversion. This difference between the even and odd cases, will lead to a difference between the maximal superintegrability manifestations between these two cases.

To discuss the consequences of the maximal superintegrability, we must enlarge the set of the initial conditions considered, allowing for a nonzero initial velocity of the light particle. Only then does the functional-independence of the new third invariant become manifest. Generally, the allowed sets of incident (in) velocities, i.e.,  the states where no collisions occurred in the past, would require positive initial velocities ordered according to
\begin{equation}
V_{\text{in}} > v_{\text{in}} > 0 \;.
\label{incident_state}
\end{equation}

Likewise, an outgoing state (out), i.e.,  a state that does not lead to any collisions in the future, is~characterized by negative final velocities ordered according to
\begin{equation}
V_{\text{out}} < v_{\text{out}} < 0\;.
\label{outgoing_state}
\end{equation}

It can be shown that the conservation of energy and the observable $J$,
\begin{eqnarray*}
T_{\text{out}}& =& T_{\text{in}}\\
J_{\text{out}}& =& J_{\text{in}}\;\,
\end{eqnarray*}
both being a function of the velocities \emph{only}, restricts the set of the allowed outgoing velocity pairs produced by the given incident pair, to one value only ({{This can be shown, in particular, by observing that }(a) the outgoing pair $(W_{\text{out}},\,w_{\text{out}})$ is an image of the incident pair, $(W_{\text{out}},\,w_{\text{out}})$, upon application of one of the elements of the group, and that (b) the condition (\ref{outgoing_state}) defines a particular \emph{chamber} of this group. However, by construction, there is only one point of an orbit of a  group per chamber}).
Notice that in this case, the outgoing velocities do not depend on the order of collisions: depending on the initial coordinates $X_0<x_0<0$, the first collision in the chain can be represented by either a ball--wall or a ball--ball collision. This independence can be regarded as a classical (as opposed to quantum) manifestation of the so-called Yang-Baxter property \cite{gaudin1983_book_english,sutherland2004_book} for the three-body system where the wall is considered a third, infinitely massive body.

\begin{figure}[h]
\centering
\includegraphics[width=0.72\columnwidth]{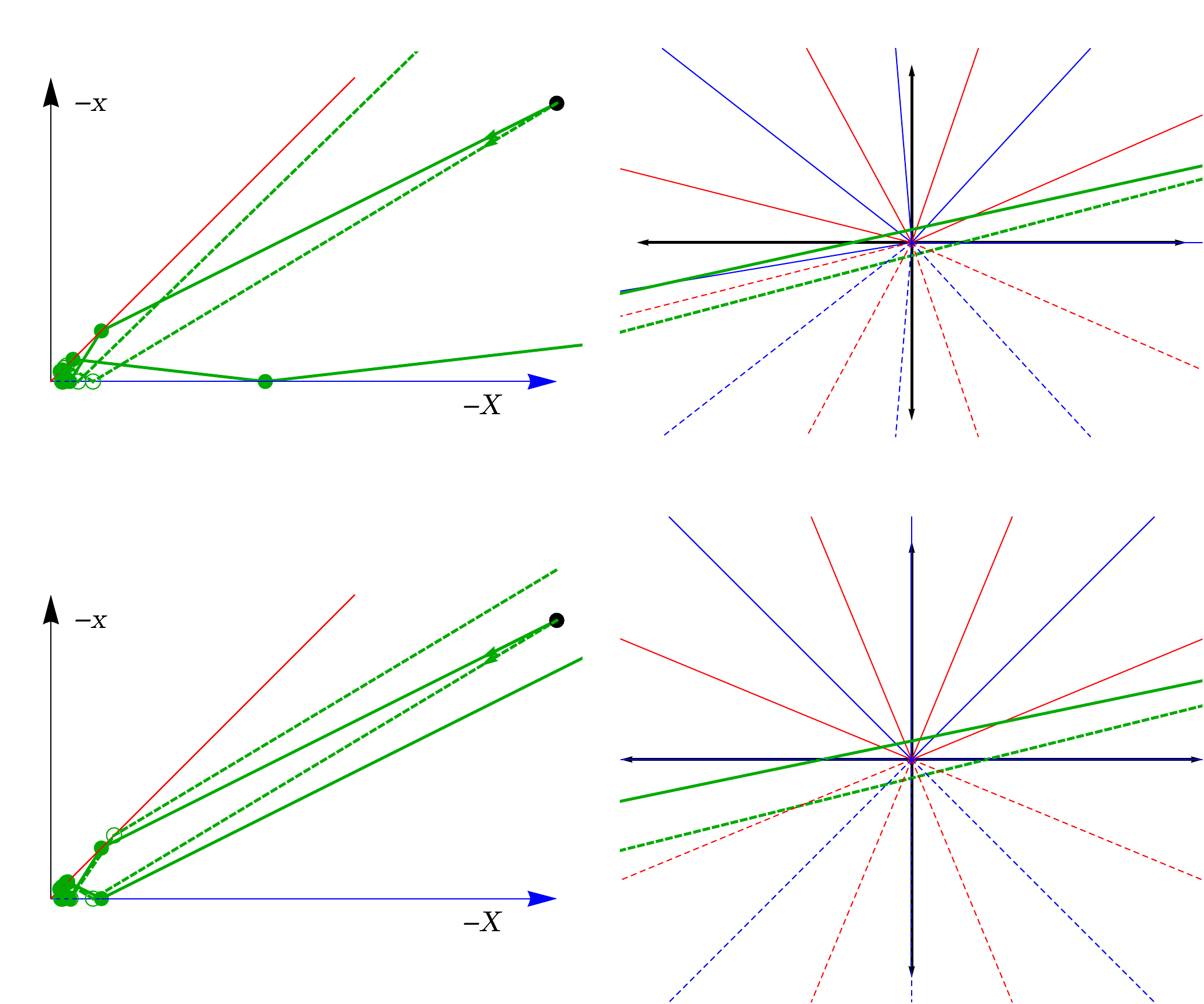}
\caption{
Each horizontal pair of figures depicts two trajectories (in green) that start at the same point in configuration space but have slightly different velocities in configuration space (left) and unfolded billiard coordinates (right). For the solid green lines, the initial  velocities determine that the first collision is a ball--ball collision (red line), and for the dashed lines the first collision is ball--wall (blue line). The top row depicts a generic mass ratio $\beta=\pi/7.6$. From the configuration space trajectory in the top left, we see that for the generic case, slightly different initial conditions can lead to very different final velocities and a different number of collisions. The figure at the top right explains this by showing that the two nearby  trajectories induce inequivalent unfoldings; the solid red and blue lines are the unfolded collision lines for the trajectory with ball--ball as first collision and the dashed red and blue lines for the trajectory with ball--wall as first collision. In this figure, the slightly different initial velocities lead to unfolded trajectories that diverge slowly and linearly in time, but the inequivalent unfoldings make their projections back to configuration space coordinates very different. In contrast, for the superintegrable case $\beta=\pi/8$ depicted in the bottom row, the solid and dashed unfoldings align at $\pi$. The unfolded trajectories with only slightly different slopes therefore give only slightly different projections at the end. Note that the reversal of initial velocities (\ref{eq:even}) is evident in the lower left figure.
}
\label{Fig:superint}
\end{figure}

In contrast to the superintegrable mass ratios, a generic mass ratio produces two different outcomes, depending on the order of collisions. Notice that two qualitatively different trajectories may even originate from two infinitely close initial conditions; see Figure~\ref{Fig:superint}. In the maximally superintegrable case of integer $q$, these two trajectories collapse to a single one-dimensional line. This~phenomenon can be regarded as an unbounded orbit analogue of the closing the orbits in the bounded case.

The actual sets of the outgoing velocities are very different in the even and in the odd cases. In~the even superintegrable case, the initial velocities are simply inverted:
\begin{equation}\label{eq:even}
q=\text{even} \to
\begin{array}{l}
\quad V_{\text{out}} = - V_{\text{in}}\\
\quad v_{\text{out}} = - v_{\text{in}}
\end{array}
\;.
\end{equation}

Indeed, since the energy and, in this case, the observable $J$ are even functions of the velocities,
the~above connection protects the conservation laws. The odd case is much more involved. One can show that
\begin{equation}
q=\text{odd} \to
\begin{array}{l}
\quad V_{\text{out}} = - \cos(\pi/q) V_{\text{in}} - \tan(\pi/q) \sin(\pi/q) v_{\text{in}}\\
\quad v_{\text{out}} = - \cos(\pi/q) ( V_{\text{in}} - v_{\text{in}})
\end{array}
\;.
\end{equation}

Remark that the case $v_{\text{in}} = V_{\text{in}}$, where $v_{\text{out}}$ vanishes, may be regarded as a generalization of a notion of a \emph{Galilean Cannon}~\cite{olshanii2016_161001060}: a system of balls that arrives at the wall with the same speed and transfers all the energy to the far-most one in the end.

\section{Physical Realizations \label{Sec:physical realizations}}
\subsection{Finite-Size Balls \label{Sec:hard rods}}

\begin{figure}[h]
\centering
\includegraphics[width=0.45\columnwidth]{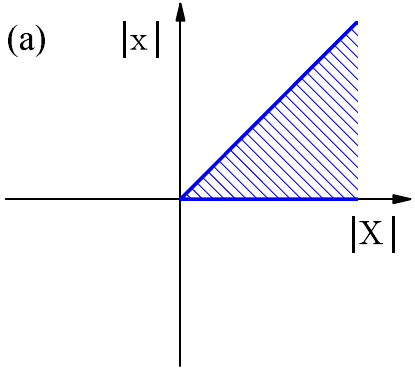}
\includegraphics[width=0.45\columnwidth]{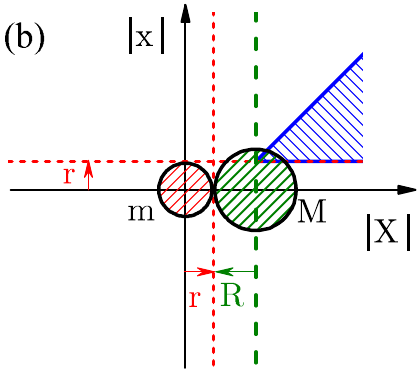}
\caption{{Configuration space} for (\textbf{a}) two point balls (\textbf{b}) balls of size $r$ and $R$.
Mapping~(\ref{Eq:excluded volume}) translates configuration space (\textbf{b}) into (\textbf{a}).
}
\label{Fig:excluded volume}
\end{figure}

The pair $(X, x)$ of positions generates a configuration point, and the set of all configuration points form the configuration space\cite{GalperinBook}. For point-size balls, it is bounded by the position of the wall, $|X|>0$ and $|x|>0$, and the condition of the impenetrability of the balls, preserving their order, $0<|x|<|X|$. More realistically, the real balls must have some finite size which we denote as $R$ and $r$ for the radii of the heavy and light balls, respectively. Still, we argue that if all collisions are elastic, the problem can be effectively reduced to the previous one of point-size balls. One might note that finite-size impenetrable balls have a smaller configuration space, schematically shown in Figure~\ref{Fig:excluded volume}, which contains an {\em excluded volume}~\cite{Girardeau60}. The configuration space of finite-size balls is $|x|>r$ and $|X|>|x|+r+R$. In other words, mapping which removes the excluded volume, 
\begin{eqnarray}
x' &=& x - r\\
X' &=& X - R - r,
\label{Eq:excluded volume}
\end{eqnarray}
reduces the problem of finite-size hard spheres to the problem of point-like objects, via a simple scaling which does not affect the balls' velocities.

As sphere is a three-dimensional object, sometimes finite width one-dimensional balls are referred to as {\em hard rods}.

\subsection{Billiard \label{Sec:billiard}}

The restricted domain of the available phase space (half of a quadrant) together with the specular reflection laws makes the system consisting of two identical balls and a wall mappable to a problem of a {\em billiard} with opening angle of $45^\circ$. In a billiard, the balls move in straight lines and collide with the boundaries (mirrors), where the incident and reflected angles are equal~\cite{KozlovTreshchevBook}. It might be shown~\cite{Galperin2003,TabachnikovBook2005} that billiard variables~(\ref{Eq:billiard variables}) change the opening angle to $\beta$ and have a special property which is that the reflections result in a straight trajectory. This unfolding creates a straight-line trajectory which intersects a certain number of lines, each of them rotated by the angle $\beta$. Each intersection corresponds to a single collision and the total number of intersections defines the total number of collisions $\Pi$. Altogether, this picture provides an intuitive visualization of the relation between $\pi$, corresponding to the angle of $180^\circ$, and the number of collisions.

\subsection{Four-Ball Chain \label{Sec:four ball chain}}

Another physical system which conceptually is related to the present system consisting of two balls and a wall, is a problem of four balls on a line. The action of the wall consists in reflecting the mass $m$ ball with the same absolute value of the velocity, $v \to -v$. The same effect can be achieved by replacing the rigid wall by another ball of mass $m$, moving with velocity $-v$. During an elastic collision, both balls will exchange their velocities. In order to make the system completely symmetric, one also has to add an additional heavy ball, resulting in $M -m -m - M$ chain. The distance between 1--2 and 3--4 balls must be the same, while 2--3 distance can be arbitrary chosen. Finally, the initial velocities should be chosen such $v_2 = v_3$ and $v_1 - v_2 = v_3 - v_4$.

\subsection{Calogero-Sutherland Particles \label{Sec:CSM}}

The form of the Hamiltonian~(\ref{eq:Calogero}) looks like a free particle in two dimensions, but note that~(\ref{eq:Calogero}) also has an interpretation as a one-dimensional particle with mass $M+m$ bouncing off a (centrifugal) barrier at $r=0$. Further, this can be mapped onto a two-particle Calogero-type model~\cite{calogero1969a,calogero1971} where $P_r$ is the relative momentum and $r$ is the relative distance. This suggests an interpretation, at least metaphorically, of the small particle  as the "force carrier" mediating an inverse-square interaction between the heavy ball and the wall.

Similar to the superintegrable (equal-mass) classical Calogero model, in which the net result of $N$ particle scattering is that the asymptotic outgoing  particle momenta are just a permutation of the incoming momenta  without any time delay~\cite{polychronakos_non-relativistic_1989}, there is also no time delay from Galperin billiards in the superintegrable cases when $\beta=\pi/q$ and $q \in 3,4,5,\ldots$. This can be seen from Equation~(\ref{eq:totalt}) and the fact that for $\beta=\pi/q$ then $\Pi\beta=\pi$ and the total time is $2|x_0|/V_0$. Including more general initial velocities $V_0$ and $v_0$ as in Section~\ref{Sec:super}, the net effect of the $\Pi$ collisions is that the two initial velocities are exactly reversed.

\section{Systematic Error \label{Sec:error}}

Any real experimental procedure should contain an error analysis. For example, the stochastic method of Buffon provides not only an approximate value of $\pi$, but also the statistical error associated with it. After $N$ trials of dropping the needle, $\pi$ is estimated as an average value while the statistical error is $\varepsilon_{stat} = \sigma / \sqrt{N-1}$, where $\sigma$ is the variance. Although in each experiment the realizations are different, the statistical error can be estimated and its value can be controllably reduced by increasing the number of trials.

\begin{figure}[h]
\centering
\includegraphics[width=0.6\columnwidth]{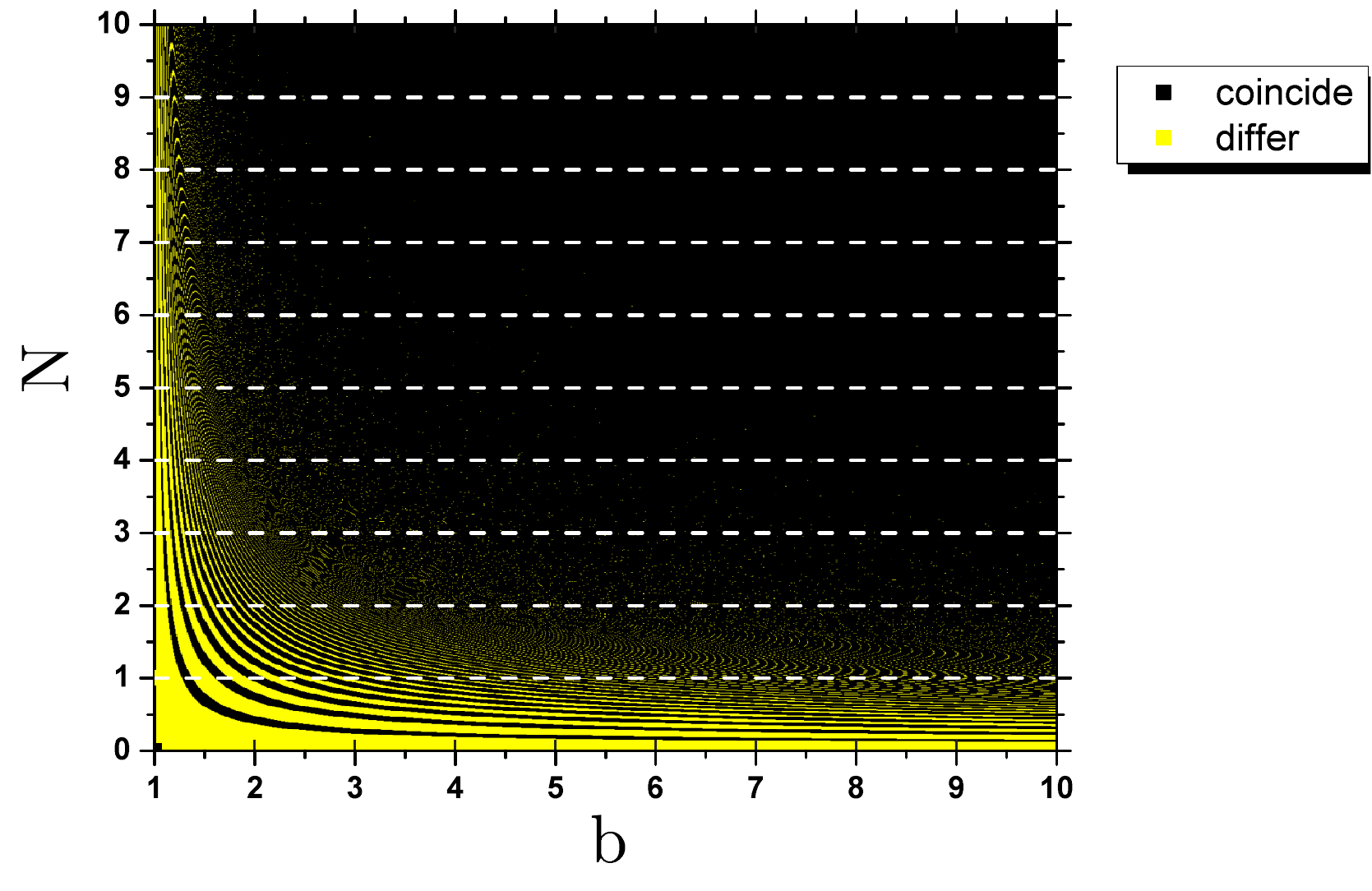}
\caption{
Difference between the exact number of collisions~(\ref{Eq:pi:acrctg}) and the approximation~(\ref{Eq:hurbilketa}), which relates it to the digits of pi,
as a function of base $b$ and mantissa $N$.
Two possible values are denoted with the dark (0) and light (1) colors.
}
\label{Fig:error}
\end{figure}

In the present study we do not report results of a real experiment, in which the number of collisions will be limited by friction, non-perfect elasticity of collisions, etc. Nevertheless, the relation~(\ref{Eq:hurbilketa}) between the number of collisions and the Galperin billiard relies on the Taylor expansion of inverse tangent function in Equation~(\ref{Eq:pi:acrctg}) and on taking its integer part, and these might induce a certain error to the final result. Accuracy of the approximations used is reported in Figure~\ref{Fig:error} as a function of the base $b$ and mantissa $N$. For completeness, here we consider $N$ not limited to integer values but as a continuous variable $N\ge 0$ and the base $b>1$. The analyzed data gives error $\varepsilon$ limited to two values $\varepsilon = 1$  (light color) and $\varepsilon = 0$ (black). It becomes evident that for large $N$ the approximate formula always works correctly, while for small $N$ there appears a complicated structure as a function of $b$. For large system base (for example, decimal $b=10$ and hexadecimal $b=16$ cases) Formula~(\ref{Eq:hurbilketa}) works correctly for any length of mantissa apart from $N=0$ case, which in any case should be treated separately due to degeneracy as will be discussed in Section~\ref{Sec:degeneracy}.

The error $\varepsilon$ is a complicated non-analytic function of $N$ and $b$, as can be perceived from Figure~\ref{Fig:error}. It turns out that for some integer bases expressions~(\ref{Eq:pi:acrctg}) and~(\ref{Eq:hurbilketa}) lead to different results. Namely, the error is $\varepsilon = 1$ for integer bases $b=6 ; 7; 14$ and $N=1$. It means that for the mentioned combinations, Galperin billiard method does not provide the digits of $\pi$ exactly, as there is an error of $\varepsilon = 1$  in the last digit. The cardinality of irrational numbers is greater than that of the integer numbers. For irrational numbers it is possible to find examples where the error is different from zero for different values of $N$ and the same value of the base $b$. Namely, $\varepsilon=1$ for $b=3.7823797$ and $N=1, 2, 3, 4$ and $6$. In general, it is clear that the closer is the base to $b=1$ the worse is the description, and for a larger number of values of $N$ Galperin billiard gives digits different from $\pi$.

We propose to treat a possible difference between~(\ref{Eq:pi:acrctg}) and~(\ref{Eq:hurbilketa}) as a {\em systematic error}, so that the final result of each ``measurement'' is $\varepsilon / b^N$ with $\varepsilon \le 1$. That is, the approximation of $\pi$ in a base $b$ can be expressed from the number of collisions $\Pi(b,N)$ as
\begin{equation}
\pi_b = \frac{\Pi(b,N)}{b^N} \pm \frac{\varepsilon}{b^N}\;.
\label{Eq:error}
\end{equation}
Such a classification is closer to a spirit of a real measurement, where different effects might contribute to the error. Another advantage of the proposed idea of introducing the concept of the systematic error, is that it solves the problem of the number of digits which are predicted correctly using Galperin billiard. It was noted by Galperin in Ref.~\cite{Galperin2003} (see also Ref.~\cite{TabachnikovBook2005}) that if there is a string of nines, that~might lead to a situation when more than one digit is different. In a similar sense, the numbers $0.999$ and $1.000$ differ by all four digits. If instead, one allows an error of $0.001$, both numbers become compatible. Indeed, from a practical point of view (suppose we calculate perimeter of a circle knowing its radius), the use of the incorrect value would lead to a relative error of 0.001 and not to completely incorrect result as all the original digits are different.

In the next sections we consider the cases of integer and non-integer bases.

\section{Integer Bases \label{Sec:integer Bases}}

Equation~(\ref{Eq:pi:acrctg}) has a profound mathematical meaning, as the number of collisions $\Pi(b,N)$ provides the first $N$ digits of the fractional part (i.e., digits beyond the radix point) of the number $\pi$ in base $b$. It might be immediately realized that as the number of collisions is obviously an integer number, its integer base representation can be chosen to be finite.

In Sections~\ref{Sec:integer Bases} and \ref{Sec:Non-integer Bases} we use number of collisions $\Pi(b, N)$, as given by Equation~(\ref{Eq:pi:acrctg}), to approximate the digits of $\pi$ for different integer bases $b$, then $(\Pi/b^N)_b$ yields the base-$b$ representation of $\pi$ with $N$ digit beyond the radix point.

\subsection{Representing a Number in Integer Bases \label{Sec:integer Representation}}

Let $b>1$ be an integer number. Any positive number $x$ has the integer expansion in base $b$, i.e., can be represented in powers of $b$ as
\begin{equation}
\label{Eq:integer}
x =(a_n a_{n-1} \ldots a_0.a_{-1}a_{-2} \ldots)_b = \sum\limits_{i=-\infty}^n a_i b^i,
\end{equation}
where $n=\integer [\log_b x]$ and $a_i=\{0,1, \ldots, b-1\}$ are the digits in the corresponding numeral system and we use form $x_b$ to denote the representation of number $x$ in base $b$. For bases with $b>10$, the symbols $A, B, \ldots$ are commonly used to denote $10, 11, \ldots$. In order to obtain the digits $a_i$, one can use the following iterative process: $a_i=\integer [ r_i/b^i ]$, $i\leq n$ with $r_n=x$ and $r_{j-1}=r_j-a_j\cdot b^j$, $j\leq n-1$. Multiplying the base-$b$ representation~(\ref{Eq:integer}) by $b^i$ shifts the radix point by $i$ digits. Thus, approximation~(\ref{Eq:hurbilketa}) gives the integer part and first $N$ digits of the fractional part of $\pi$ in base $b$.

The most frequently used integer systems are decimal ($b=10$) and binary $b=2$ systems. Occasionally, also ternary $b=3$, octal ($b=8$), hexadecimal ($b=16$) and others systems are used. Importantly, for integer bases, finite representations are unique, while infinite representations might be not unique. For example, the finite number $1_{10}$ in the decimal base can be written as $1.000(0)_{10}=0.999(9)_{10}$.

\subsection{Degenerate Case of Equal Masses and Submultiple Angles \label{Sec:degeneracy}}

Before considering in detail the representation in bases $b=10; 2; 3$ reported in Tables~\ref{table:b=10}--\ref{table:b=3}, we~note that $N=0$ case is universal as the mass ratio $M/m = b^N = 1$ does not depend on the base $b$. In other words, the digit in front of the radix point always correspond to the same number. The~Equation~(\ref{Eq:pi:acrctg}) formally gives 4 collisions, which is different from the physically correct number of 3 collisions. The~reason for such a difference comes from a degeneracy between the third and fourth collision. While for $N>0$, the direction of the light ball is always inverted in the last two collisions ($\phi\to -\phi$), for $N=0$ the light ball completely stops exactly at the third collision. In physical sense there is no difference between $v_3 = -0$ and $v_4 = +0$ velocities.

Thus, Equations~(\ref{Eq:pi:acrctg}) and (\ref{Eq:hurbilketa}) are applicable only for $N\geq1$ while $N=0$ is a special case and it should be treated separately.

The analogous result takes place in the case of the angle $\phi$ being submultiple of $\pi$, i.e., when~the ratio $\pi/\phi$ is an integer number. The number of collisions is not given correctly by~(\ref{Eq:pi:acrctg}) as the last collision is degenerate as well.

\subsection{Decimal Base \label{Sec:Decimal Base}}

For the most natural case of the decimal base system, $b=10$, the number of collisions $\Pi(10, N)$ is given in Table~\ref{table:b=10}. It is easy to follow, how Galperin billiard generates digit of $\pi$. For $N=0$, Equation~(\ref{Eq:pi:acrctg}) results in the first digit of $\pi$ approximated by 4, while due to degeneracy discussed in Section~\ref{Sec:degeneracy}, physically there are 3 collisions. For $N=1$, there are $31$ collisions, resulting in expression with 1 digit after the radix point, $3.1$. From $N=2$, the number of collisions in 314 giving the number $\pi$ with 2 digits after the radix point. One can see that the billiard method correctly approximates the number $\pi$ as 3 plus $N$ more digits in the decimal base.

Conceptually, one might ask if there is a difference between the number of collisions Equation~(\ref{Eq:pi:acrctg}) which depend on $\arctan(b^{-N})$ rather than $b^{-N}$, as in Equation~(\ref{Eq:hurbilketa}). It turns out that the base $b=10$ is large enough (see Figure~\ref{Fig:error}) so that there is no any difference in the integer part of the expansion.

\begin{table}[h]
\centering
\caption{Number of collisions $\Pi(10,N)$ given by Equation~(\ref{Eq:pi:acrctg}) for the decimal base, $b=10$.
The first column reports the value of mantissa $N$.
The second column is the resulting number of collisions in the decimal base.
The third column is the number $\pi$ with $N$ digits in the fractional part in the decimal representation.
The fourth column gives the systematic error according to Equation~(\ref{Eq:error}). The case where approximation~(\ref{Eq:hurbilketa}) fails as compared to~(\ref{Eq:pi:acrctg}) is highlighted by red. 
The blue digit is incorrectly predicted by the Galperin billiard. 
}
\label{table:b=10}
\begin{tabular}{clll}
\toprule
\textbf{\emph{N}}  & \boldmath{$\Pi(10,N)_{10}$}       & \boldmath{$(\Pi(10,N)/10^N)_{10}$  } & \boldmath{$(1/10^N)_{10}$}
\\ \midrule
0  & \color{red}{4}           & \color{blue}4            & 1 \\

1  & 31          & 3.1          & 0.1\\
2  & 314         & 3.14         & 0.01\\
3  & 3141        & 3.141        & 0.001\\
4  & 31415       & 3.1415       & 0.0001\\
5  & 314159      & 3.14159      & 0.00001\\
6  & 3141592     & 3.141592     & 0.000001\\
7  & 31415926    & 3.1415926    & 0.0000001\\
8  & 314159265   & 3.14159265   & 0.00000001\\
9  & 3141592653  & 3.141592653  & 0.000000001\\
10 & 31415926535 & 3.1415926535 & 0.0000000001\\
\bottomrule
\end{tabular}
\end{table}
\unskip
\begin{table}[h!]
\centering
\caption{
Number of collisions $\Pi$ given by Equation~(\ref{Eq:pi:acrctg}) for the binary base, $b=2$.
The first column reports the value of mantissa $N$.
The second column is the resulting number of collisions in the decimal base.
The third column is the number of collisions written in the binary representation.
The fourth column is the binary representation of the number $\pi$ with $N$ digits in the fractional part.
The fifth column gives the systematic error according to Equation~(\ref{Eq:error}). 
The case where approximation~(\ref{Eq:hurbilketa}) fails as compared to~(\ref{Eq:pi:acrctg}) is highlighted by red. 
The blue digits are not correctly predicted by the Galperin~billiard.
}
\label{table:b=2}
{
\begin{tabular}{lllll}
\toprule
\textbf{\emph{N}} & \boldmath{$\Pi(2,N)_{10}$} & \boldmath{$\Pi(2,N)_{2}$}  & \boldmath{$(\Pi(2,N)/2^N)_{2}$}
& \boldmath{$(1/2^N)_{2}$}

\\ \midrule
0 & \color{red}4     & 1\color{blue}00        & 1\color{blue}00            & 1\\
1 & 6     & 110        & ~11.0           & 0.1\\
2 & 12    & 1100       & ~11.00          & 0.01\\
3 & 25    & 11001      & ~11.001         & 0.001\\
4 & 50    & 110010     & ~11.0010        & 0.0001\\
5 & 100   & 1100100    & ~11.00100       & 0.00001\\
6 & 201   & 11001001   & ~11.001001      & 0.000001\\
7 & 402   & 110010010  & ~11.0010010     & 0.0000001\\
8 & 804   & 1100100100 & ~11.00100100    & 0.00000001\\
9 & 1608  & 11001001000& ~11.001001000   & 0.000000001\\
10 & 3216 & 110010010000& ~11.0010010000 & 0.0000000001\\
\bottomrule
\end{tabular}
}
\end{table}
\unskip
\begin{table}[h!]
\centering
\caption{
Number of collisions $\Pi$ given by Equation~(\ref{Eq:pi:acrctg}) for the ternary base, $b=3$.
The first column reports the value of mantissa $N$.
The second column is the resulting number of collisions in the decimal base.
The third column is the number of collisions written in the binary representation.
The fourth column is the ternary representation of the number $\pi$ with $N$ digits in the fractional part.
The fifth column gives the systematic error according to Equation~(\ref{Eq:error}). 
The case where approximation~(\ref{Eq:hurbilketa}) fails as compared to~(\ref{Eq:pi:acrctg}) is highlighted by red. 
The blue digits are incorrectly predicted by the Galperin~billiard.
}
\label{table:b=3}{
\begin{tabular}{lllll}\toprule
\textbf{\emph{N}} & \boldmath{$\Pi(3,N)_{10}$} & \boldmath{$\Pi(3,N)_{3}$}  & \boldmath{$(\Pi(3,N)/3^N)_{3}$}
& \boldmath{$(1/3^N)_{3}$}\\ \midrule
0 &\color{red} 4 & 1\color{blue}1 & 1\color{blue}1 & 1\\
1 & 9 & 100 & 10.0 & 0.1\\
2 & 28 & 1001  & 10.01 & 0.01\\
3 & 84 & 10010 & 10.010 & 0.001\\
4 & 254 & 100102 & 10.0102 & 0.0001\\
5 & 763 & 1001021  & 10.01021 & 0.00001\\
6 & 2290 & 10010211 & 10.010211 & 0.000001\\
7 & 6870 & 100102110  & 10.0102110 & 0.0000001\\
8 & 20611 & 1001021101 & 10.01021101 & 0.00000001\\
9 & 61835 & 10010211012  & 10.010211012 & 0.000000001\\
10 & 185507 & 100102110122 & 10.0102110122 & 0.0000000001\\\bottomrule
\end{tabular}
}
\end{table}

\subsection{Binary and Ternary Bases \label{Sec:binary Base}}

Other important examples of number systems include the binary ($b=2$) and ternary ($b=3$) base systems. The binary system lies in the core of modern computers which operate with {\em bits} $0, 1$. Interestingly, base-$3$ computer named Setun was built 1958 under leadership of mathematician Sergei Sobolev and operated with {\em trits}, $0,1,2$.

Table~\ref{table:b=2} reports the number of collisions $\Pi(2,N)$ obtained for $b=2$ base. By expressing the number of collisions in binary base using zeros and ones, one obtains the representation of the number $\pi$ in binary base. In the ternary base, the number of collisions are written using the three allowed digits, $0,1,2$, see Table~\ref{table:b=3}.

\subsection{Best Bases for a Possible Experiment \label{Sec:experiment}}

As concerning the effects of the friction and other sources of energy dissipation, it is easier to perform experiments for small base $b$.
While for $N=0$ (the mass ratio is $M/m = 1$ independently of $b$) there are 3 collisions which can be easily observed with identical balls, for larger $N$ the number of collisions grows exponentially fast.
The decimal system has a rather ``large'' base $b=10$ which already for $N=1$ results in 31 collisions and $N=2$ even in 314 collisions.
It might be notoriously hard to create a clean system in which such a large number of collisions can be reliably measured.

For the binary base $b=2$ and $N=1$ the number of collisions to be observed is much smaller, $3; 6; 12; 25; \ldots$, see Table~\ref{table:b=2}, making such system more suitable for an experimental observation.

\section{Non-Integer Bases \label{Sec:Non-integer Bases}}

As anticipated above, the Galperin billiard method should provide digits of $\pi$ in an arbitrary base $b$, even for a non-integer one. In this Section we consider a number of examples.

\subsection{Representing a Number in a Non-Integer Base \label{Sec:representation}}

For a non-integer base $b>1$, any positive number $x$ can be written in the base-$b$ representation similarly to Equation~(\ref{Eq:integer})
where digits $a_i$ can take only non-negative integer values smaller than non-integer base, $a_i<\lceil b \rceil$ ($\lceil x \rceil$ stands for the least integer which is greater than or equal to $x$).

Unlike the integer bases, for a non-integer base $b$, even finite fractions might have different $b$ representations. For example, in the golden mean $\varphi\approx 1.61803$ base, due to the quadratic equality $\varphi^2=\varphi+1$, one has $100_\varphi = 11_\varphi$. With Equation~(\ref{Eq:integer}) we can find at least one representation for $x$. Moreover, the number of collisions $\Pi(b,N)$ which is obviously an integer number and  always written with a finite representation in any integer base system, in a non-integer base $b$ it is a common situation that an integer number needs an infinite representation which corresponds to the number. An example of a rational base is displayed in Table~\ref{table_1p5} for $b=3/2$. Thus, $3/2$ representation is $\pi_{3/2}=(3/2)^2+(3/2)^{-1} +(3/2)^{-4} + (3/2)^{-9} \ldots$. Note a peculiarity of this base is that $b$ is so small that the expansion $\arctan (b^{-N}) \approx b^{-N}$ used to derive Equation~(\ref{Eq:hurbilketa}) is not precise enough. As a result, the number of collisions $\Pi(b, 1)$ in the Galperin billiard, given by Equation~(\ref{Eq:pi:acrctg}), for several values of $N$ (highlighted by red in Table~\ref{table_1p5}) does not coincide with expression~(\ref{Eq:hurbilketa}) which is used to transcribe $N$ digits of $\pi$ in the base $b$. The resulting possible difference in the last digits (highlighted by blue) can be interpreted as the systematic error $\varepsilon=1$ in the spirit of Section~\ref{Sec:error}.
\begin{table}[h!]
\centering
\caption{Number of collisions $\Pi(3/2,N)$ given by (\ref{Eq:pi:acrctg}) and approximation of $\pi$ for $b=3/2$.
The first column is $N$.
The second column is $\Pi(3/2,N)$ in the decimal base.
The third column is the integer part of the number of collisions $\Pi(3/2,N)$ written in the base $3/2$.
The fourth column is the number $\pi$ with $N$ digits in the fractional part in the base $3/2$.
The fifth column gives the systematic error according to Equation~(\ref{Eq:error}). 
The cases where approximation~(\ref{Eq:hurbilketa}) fails as compared to~(\ref{Eq:pi:acrctg}) are highlighted by red. 
The blue bold digits are predicted incorrectly by the Galperin billiard.
}
\label{table_1p5}{
\begin{tabular}{lllll}\toprule
\textbf{\emph{N}} & \boldmath{$\Pi(3/2,N)_{10}$} & \boldmath{$\Pi(3/2,N)_{3/2}$}  & \boldmath{$(\Pi(3/2,N)/(3/2)^N)_{3/2}$} & \boldmath{$(1/(3/2)^N)_{3/2}$} \\ \midrule
0 &\bf\color{red} 4 & 100\bf\color{blue}0.  & 100\bf\color{blue}0. & 1\\
1 &\bf\color{red} 5 & 10\bf\color{blue}10.  & ~10\bf\color{blue}1.0  & 0.1\\
2 & 7 & 10010. & ~100.10  &  0.01\\
3 & 10 & 100100. & ~100.100  & 0.001\\
4 &\bf\color{red} 16 & 1001001. & ~100.1001  & 0.0001\\
5 & 23 & 100100\bf\color{blue}00. & ~100.100\bf\color{blue}00  & 0.00001\\
6 & 35 & 100100\bf\color{blue}010. & ~100.100\bf\color{blue}010 & 0.000001\\
7 & 53 & 100100\bf\color{blue}0100. & ~100.100\bf\color{blue}0100  & 0.0000001\\
8 & 80 & 10010010000. & ~100.10010000  & 0.00000001\\
9 & 120 & 100100100000. & ~100.100100000  & 0.000000001\\
10 & 181 & 1001001000001. & ~100.1001000001  & 0.0000000001\\\bottomrule
\end{tabular}
}
\end{table}

\subsection{Number Systems with Irrational Bases}\label{Sec:irrational base}

Some notable examples of non-integer bases include the fundamental cases of the bases with the quadratic number $b=\sqrt{3}=1.732\cdots$, the golden mean $b=\varphi=1.618\cdots$, the natural logarithm $b=e=2.718\cdots$ and a curious situation when the number $\pi$ is used itself as a base, $b=\pi$.

Table~\ref{table_sqrt3} contains the results of $\sqrt{3}$ representations of $pi$ number. Note that this includes a superintegrable value of the ratio of masses, see Section~\ref{Sec:super}. Namely, for $N=1$, the angle $\beta=\pi/6$ in~(\ref{Eq:pi:acrctg}) and $M/m=b^{2N}=3$. In this case, the values $\Pi(\sqrt{3},1)$, given by (\ref{Eq:pi:acrctg}) and (\ref{Eq:hurbilketa}), differ by 1, since~in~(\ref{Eq:pi:acrctg}) under the integer part there is an exact integer number and, therefore, an approximation of the $\arctan$ function will inevitably give an error. As for $N=5$,  the difference between the values obtained by (\ref{Eq:pi:acrctg}) and (\ref{Eq:hurbilketa}) is also 1 due to the same nature of systematic errors as for $N=1$ and $N=4$ of $b=3/2$ in Table~\ref{table_1p5}.

In Table~\ref{table_phi} we give two different representations for $\pi_\varphi$, since its integer part $100_\varphi=11_\varphi$. The~fourth and fifth columns in Table report the $\varphi$-representation with the integer part $100_\varphi$ and $11_\varphi$, respectively. The allowed digits for both representations are $0$ and $1$. Table~\ref{table_e} reports the resulting number of collisions $\Pi(e, N)$ in base $e$ with the allowed digits $0, 1$ and $2$. We get the representation $\pi=(10.1010020200\ldots)_e=e+e^{-1}+e^{-3}+2 e^{-6} + 2 e^{-8}+\ldots$. One can see the influence of the error of the computation by the Galpelin billiard which in the case of base $b=e$ is $1/e^N=0.00\ldots01_e$. Due to this error, the last digit may be incorrectly predicted by the method. Especially, when the last digit is the maximum allowed ($2$ for $b=e$, see the cases $N=4$ and $N=9$ in Table~\ref{table_e}), then the two last digits may be incorrect.
\begin{table}[h!]
\centering
\caption{Number of collisions $\Pi(\sqrt{3},N)$ given by (\ref{Eq:pi:acrctg}) and (\ref{Eq:hurbilketa}) and approximation of $\pi$ for $b=\sqrt{3}$.
The first column is $N$.
The second column is $\Pi(\sqrt{3},N)$ calculated by (\ref{Eq:pi:acrctg}) in the decimal base.
The third column is the integer part of the number of collisions $\Pi(\sqrt{3},N)$ given by (\ref{Eq:pi:acrctg}) and written in the base $\sqrt{3}$.
The fourth column is the number $\pi$ with $N$ digits in the fractional part in the base $\sqrt{3}$ calculated by (\ref{Eq:hurbilketa}) .
The fifth column is $\Pi(\sqrt{3},N)$ given by (\ref{Eq:hurbilketa}) in the decimal base. The sixth column is the integer part of the number of collisions $\Pi(\sqrt{3},N)$ given by (\ref{Eq:hurbilketa}) and written in the base $\sqrt{3}$. The seventh column is  the number $\pi$ calculated by (\ref{Eq:hurbilketa}) with $N$ digits in the fractional part in the base $\sqrt{3}$. The eighth column  gives the systematic error according to Eq.~(\ref{Eq:error}). 
The cases where approximation~(\ref{Eq:hurbilketa}) fails as compared to~(\ref{Eq:pi:acrctg}) are highlighted by red. 
The blue color stands for the digits which are incorrectly predicted by the Galperin billiard.
}
\label{table_sqrt3}{
\begin{tabular}{l|l|l|l|l|l|l|l}
\toprule
\rotatebox{-90}{$N$}
& 
\rotatebox{-90}{$\Pi_{10}$}
& 
~~~~~~~~~\rotatebox{-90}{$\Pi_{\sqrt{3}}$}
& 
~~~~~~\rotatebox{-90}{$\left(\frac{\Pi}{\sqrt{3}^N}\right)_{\sqrt{3}}$}
& 
\rotatebox{-90}{$\Pi(\sqrt{3}, N)_{10}$}
& 
~~~~~~~~~~\rotatebox{-90}{$\Pi(\sqrt{3},N)_{\sqrt{3}}$}
&
~~~~~~~\rotatebox{-90}{$\left(\frac{(\Pi(\sqrt{3}, N)}{\sqrt{3}^N}\right)_{\sqrt{3}}$~~~~~~~~}
&
~~~\rotatebox{-90}{$\left(\frac{1}{\sqrt{3}^N}\right)_{\sqrt{3}}$}
\\ \hline
0 &\bf\color{red} 4 & 10\bf\color{blue}1. & 10\bf\color{blue}1. &\bf\color{red} 3 & 100. & 100. & 1\\
1 & \bf\color{red}6 & 1000. & 100.0 &\bf\color{red} 5 & 1\bf\color{blue}10. & 1\bf\color{blue}1.0 & 0.1\\
2 & 9 & 10000. & 100.00 & 9 & 10000. & 100.00 & 0.01\\
3 & 16 & 100000. & 100.000 & 16 & 100000.  & 100.000 & 0.001\\
4 & 28 & 1000001. & 100.0001  & 28 & 1000001.  & 100.0001 & 0.0001 \\ 
5 &\bf\color{red} 49 & 10000010. & 100.00010 &\bf\color{red} 48 & 100000\bf\color{blue}01. & 100.000\bf\color{blue}01 & 0.00001\\
6 & 84 & 100000100. & 100.000100 & 84 & 100000100. & 100.000100 & 0.000001\\
7 & 146 & 1000001000. & 100.0001000  & 146 & 1000001000. & 100.0001000 & 0.0000001\\
8 & 254 & 10000010010. & 100.00010010  & 254 & 10000010010. & 100.00010010 & 0.00000001\\
9 & 440 & 100000100100. & 100.000100100  & 440 & 100000100100. & 100.000100100 & 0.000000001 \\
10\!\!& 763 & 1000001001010.\!\!& 100.0001001010\!\!& 763 & 1000001001010.\!\!& 100.0001001010\!\!& 0.0000000001\!\!\\
\end{tabular}
	}
\end{table}

In Table~\ref{table_pi} we show approximations of $\pi$ in the base $b=\pi$
The allowed digits in this base are $0, 1, 2$ and $3$.
The Galperin billiard does not provide an integer-number representation for the number $\pi$ even in this case, as instead of the ``natural'' possibility $\pi = 10_\pi$ one obtains an infinitely long representation
$$\pi = (3.0110211100\ldots)_{\pi}=3+\pi^{-2}+\pi^{-3}+2\pi^{-5}+\pi^{-6} + \pi^{-7} + \pi^{-8} +\ldots.$$

This non-unique representation is similar to the infinite representation $0.999(9)\ldots$ of $1$ in the decimal system. We note the difference between the values obtained by (\ref{Eq:pi:acrctg}) and (\ref{Eq:hurbilketa}) is also 1 for $N=1$ due to a systematic error.
\begin{table}[h!]
\centering
\caption{Number of collisions $\Pi(\varphi,N)$ given by (\ref{Eq:pi:acrctg}) for $b=\varphi$.
The first column is $N$.
The second column is $\Pi(\varphi,N)$  in the decimal base.
The third column is the integer part of $\Pi(\varphi,N)$ written in the base $\pi$.
The fourth column is  the number $\pi$ with $N$ digits in the fractional part (Type I) in the base $\pi$.
The fifth column is the number $\pi$ with $N$ digits in the fractional part (Type II) in the base $\pi$.
The~sixth column gives the systematic error according to Equation~(\ref{Eq:error}). 
The case where approximation~(\ref{Eq:hurbilketa}) fails as compared to~(\ref{Eq:pi:acrctg}) is highlighted by red. 
The~blue bold digit is predicted incorrectly by the Galperin~billiard.
}
\label{table_phi}{
\begin{tabular}{llllll}\toprule
\emph{\textbf{N}} & \boldmath{$\Pi(\varphi, N)_{10}$} & \boldmath{$\Pi(\varphi, N)_{\varphi}$}  & \boldmath{$(\Pi(\varphi, N)/\varphi^N)_{\varphi}$} \textbf{(I)} & \boldmath{$(\Pi(\varphi, N)/\varphi^N)_{\varphi}$} \textbf{(II)} & \boldmath{$(1/\varphi^N)_{\varphi}$}  \\ \midrule
0 &\bf\color{red} 4 & 10\bf\color{blue}1. & 10\bf\color{blue}1. & {\bf\color{blue}10}1.  & 1\\
1 & 5 & 1000. & 100.0 & 11.0  & 0.1\\
2 & 8 & 10001.  & 100.01 & 11.01 & 0.01\\
3 & 13 & 100010.  & 100.010 & 11.010 & 0.001\\
4 & 21 & 1000100.  & 100.0100 & 11.0100 & 0.0001\\
5 & 34 & 10001000.  & 100.01001 & 11.01001 & 0.00001\\
6 & 56 & 100010010.  & 100.010010 & 11.010010 & 0.000001\\
7 & 91 & 1000100101.  & 100.0100101 & 11.0100101 & 0.0000001\\
8 & 147 & 10001001010.  & 100.01001010 & 11.01001010 & 0.00000001\\
9 & 238 & 100010010100.  & 100.010010101 & 11.010010101 & 0.000000001\\
10 & 386 & 1000100101010.  & 100.0100101010 & 11.0100101010 & 0.0000000001\\\bottomrule
\end{tabular}
}
\end{table}
\unskip
\begin{table}[h!]
\centering
\caption{Number of collisions $\Pi(e,N)$ given by (\ref{Eq:pi:acrctg}) and approximation of $\pi$ for $b=e$.
The first column is $N$.
The second column is $\Pi(e,N)$ in the decimal base.
The third column is the integer part of the number of collisions $\Pi(e,N)$ written in the base $e$.
The fourth column is the number $\pi$ with $N$ digits in the fractional part in the base $e$.
The fifth column gives the systematic error according to Equation~(\ref{Eq:error}). 
The case where approximation~(\ref{Eq:hurbilketa}) fails as compared to~(\ref{Eq:pi:acrctg}) is highlighted by red. 
The blue bold digits are incorrectly predicted  by the Galperin billiard due to a systematic error.
}
\label{table_e}{
\begin{tabular}{lllll}\toprule
\emph{\textbf{N}} & \boldmath{$\Pi(e,N)_{10}$} & \boldmath{$\Pi(e,N)_{e}$}  & \boldmath{$(\Pi(e,N)/e^N)_{e}$} & \boldmath{$(1/e^N)_{e}$} \\ \midrule
0 & \color{red}4      &1\color{blue}1. &1\color{blue}1. & 1\\
1 & 8      & 100. & 10.0 & 0.1\\
2 & 23     & 1010. & 10.10 & 0.01\\
3 & 63     & 10101. & 10.101 & 0.001\\
4 & 171    & 1010{\color{blue} 02}. & 10.10{\color{blue}  02} & 0.0001\\
5 & 466    & 1010100. & 10.10100 & 0.00001\\
6 & 1267   & 1010100{\color{blue}1}. & 10.10100{\color{blue}1} & 0.000001\\
7 & 3445   & 101010020. & 10.1010020 & 0.0000001\\
8 & 9364   & 1010100201. & 10.10100201 & 0.00000001\\
9 & 25456  & 101010020{\color{blue}  12}. & 10.1010020{\color{blue}  12} & 0.000000001\\
10 & 69198 & 101010020200. & 10.1010020200 & 0.0000000001\\\bottomrule
\end{tabular}
}
\end{table}

The accuracy of the approximation to $\pi$ obtained by Galperin's billiard can be estimated by Equation~(\ref{Eq:error}) which we treat as systematic error of the method. It can be noted that the wrong digits might appear due to failure of approximation~(\ref{Eq:hurbilketa}) as compared to~(\ref{Eq:pi:acrctg}) which is highlighted by red in Tables~\ref{table:b=10}--\ref{table_pi}, either by differences coming due representation of an integer number (i.e., the~number of collisions) in a non-integer base (see Tables~\ref{table_1p5},\ref{table_sqrt3},\ref{table_e}).
\begin{table}[h!]
\centering
\caption{Number of collisions $\Pi(\pi,N)$ given by (\ref{Eq:pi:acrctg}) for $b=\pi$.
The first column is $N$.
The second column is the number of collisions in the decimal base.
The third column is the integer part of the number of collisions $\Pi(\pi,N)$ written in the base $\pi$.
The fourth column is the number $\pi$ with $N$ digits in the fractional part in the base $\pi$.
The fifth column gives the systematic error according to Equation~(\ref{Eq:error}). The case $N=1$ is emphasized since there is the difference 1 between $\Pi(\pi,1)$ by~(\ref{Eq:pi:acrctg}) and the approximation~(\ref{Eq:hurbilketa}). 
The cases where approximation~(\ref{Eq:hurbilketa}) fails as compared to~(\ref{Eq:pi:acrctg}) are highlighted by red. 
The blue bold digits are not correctly predicted  by the Galperin billiard.
}
\label{table_pi}{
\begin{tabular}{lllll}\toprule
\emph{\textbf{N}} & \boldmath{$\Pi(\pi,N)_{10}$} & \boldmath{$\Pi(\pi,N)_{\pi}$} & \boldmath{$(\Pi(\pi, N)/\pi^N)_{\pi}$} & \boldmath{$(1/\pi^N)_{\pi}$} \\ \midrule
0 &\bf\color{red} 4       & \color{blue}10. & \color{blue}10. & 1\\
1 &\bf\color{red} 10      & \color{blue}100. &\color{blue} 10.0 & 0.1\\
2 & 31      & 301. & 3.01 & 0.01\\
3 & 97      & 3010. & 3.010 & 0.001\\
4 & 306     & 30110. & 3.0110 & 0.0001\\
5 & 961     & 301102. & 3.01102 & 0.00001\\
6 & 3020    & 3011021. & 3.011021 & 0.000001\\
7 & 9488    & 30110210. & 3.0110210 & 0.000001\\
8 & 29809   & 301102110. & 3.01102110 & 0.0000001\\
9 & 93648   & 3011021110. & 3.011021110 & 0.00000001\\
10 & 294204 & 30110211100. & 3.0110211100 & 0.000000001\\\bottomrule
\end{tabular}
}
\end{table}

\vspace*{2cm}
\section{Conclusions \label{Sec:conclusions}}

To summarize, we have studied how the digits of the number $\pi$ are generated in a simple classical three-body system consisting of one heavy ball, one light ball and a wall (Galperin's billiard method).
We obtain for the first time, to our best knowledge, the complete explicit solution for the balls' positions and velocities as a function of the collision number and time. This is achieved by moving to billiard coordinates and unfolding the trajectory. In this representation, collisions are reflections and the motion looks almost like a free particle moving in two dimensions if the singular impulse at reflections is ignored. The square of the angular momentum about the three-body coincidence of the two balls with the wall is an invariant integral of motion. This quantity explains why the invariant first identified the adiabatic approximation works not only close to the return point but also for any ball--wall collision, even far from the wall.

This free-particle form of the Hamiltonian~(\ref{eq:Calogero}) which looks like the Calogero model also explains many of the ``round'' properties of Galperin billiards, such as the portraits of the system in $(P,1/X)$ and $(V,1/X)$ coordinates have a shape close to a circular one. Another circle appears in $(V,v)$ coordinates and corresponds to the energy conservation law.
Instead a hyperbolic shape appears in the $(X,t)$ plane. A third invariant is also revealed for certain mass ratios that makes the system superintegrable and removes the dependence of the number of collisions on the initial conditions for generalized scenarios.

A recent article establishing an isomorphism between the dynamics in Galperin billiards and Grover's algorithm for quantum database searches~\cite{brown} inspires consideration of the quantum version.  Curiously, since the Galperin model effectively becomes a classical simulator for a quantum algorithm, its quantum realization would lead to a \emph{second quantization of the Grover process}, where the coefficients 
 of Grover's wavefunction will be promoted 
 to quantum observables. Note that unlike in the standard second quantization of field theory, components of the first-quantized wave 
 function of the space Grover's machine lives in become \emph{Hermitian}
 operators, ball velocities in our case. On the negative side, while 
 the map presented in~\cite{brown} is indeed a one-to-one correspondence 
 between the two protocols, as~far as the degrees 
 of freedom are concerned,
 the only isomorphism it establishes is a map between Galperin's velocity space and a two-dimensional Hilbert space Grover's algorithm constrains the database to: the actual Hilbert space the database lives in can be arbitrary large. Nonetheless, potential advantages of a second quantization on this reduced Hilbert space are worth exploring.

 Another motivation for considering the quantum version of this system is the daunting  experimental challenge of realizing Galperin billiards with macroscopic objects.  Quantum realizations of effectively one-dimensional mixed-mass systems like Galperin billiards can be experimentally realized with ultracold atoms in few-body systems~\cite{few} or as bi-solitons in ultracold atomic gases, via a scheme described in~\cite{olshanii2016_161001060}. There, a bi-soliton of a desired mass ratio is created using a coupling constant quench~\cite{satsuma}; one of the two solitons is subsequently transferred to a different internal atomic state that leads to a repulsion between the solitons. The wall is generated using a light sheet. High degree of macroscopic quantum coherence will be guaranteed \cite{yurovsky2017_220401} by the higher conservation laws operating in the one-dimensional bosonic systems. Bi-solitons have been recently created experimentally \cite{haller,hulet_private}. 
 
 The integrals of motion, including the third superintegral, should carry over without modification into quantum observables.  In fact, mixed-mass superintegrability with hard-core interactions has been previously identified for quantum billiards in free space~\cite{olshanii2015_105005} and harmonic traps~\cite{harshman2017}. 
 
Examples of integer bases $b$, including decimal, binary and ternary, are considered. We argue that smaller bases (for example, $b=2;3$) are the easiest to be realized in an experiment and show how the Galperin billiard can be generalized to finite-size balls (hard rods). We show that the dependence of a possible error in the last digit as a function of $b$ and $N$ has a complicated form, with the error disappearing in the limit of $b^N\to\infty$. We propose to treat the possible error in the last digit as a systematic error. In particular this resolves the problem of the correct number of obtained digits. Finally, we consider non-integer bases, including expressing $\pi$ in the base $\pi$ or the base of the golden ratio $\varphi$. These reveal the curious limitations of numeral representations of irrational numbers in irrational bases, and make Galperin's $\pi$-calculating machine all the more remarkable.

\vspace{6pt} 

\section*{Acknowledgment}
{The research leading to these results received funding from the MICINN (Spain) Grant No.  FIS2017-84114-C2-1-P.
M.G. has been partially supported by Juan de la Cierva-Formaci\'on FJCI-2014-21229 and Juan de la Cierva-Incorporaci\'on IJCI-2016-29071 fellowships, the Spanish grants 	MICIIN/FEDER MTM2015-65715-P, MTM2016-80117-P (MINECO/FEDER, UE), PGC2018-098676-B-I00 (AEI/FEDER/UE) and the Catalan grant 2017SGR1374. M.O. acknowledges financial support from the National Science Foundation grants PHY-1607221 and PHY-1912542 and the US-Israel Binational Science Foundation grant 2015616.}

\appendix
\section{ Solving the Equations of Motion \label{Sec:equations of motion}}
\unskip

In Section~\ref{Sec:number of collisions}, it was shown that the total number of collisions can be explicitly obtained from the conservation laws resulting in Equation~(\ref{Eq:hurbilketa}), and it does not depend on the exact initial position of the lighter ball $x_0$ or the initial velocity of the incident ball $V_0$. Here we outline how the trajectory of the balls can be obtained.

First, one has to integrate the equations of motion, as for example in \cite{Weidman2013}. Let $X_n$ and $x_n$ be the coordinates and velocities of the heavy and light balls at the time $t_n$ of the $n$-th collision, respectively. As the velocities change only at contact, the balls move with constant velocities $V_n$ and $v_n$ between the $n$-th and $(n+1)$-th collisions. The velocities change according to the rules (the energy and momentum conservation laws) provided in Section~\ref{Sec:number of collisions}.
All odd collisions,  $n=2k+1$, correspond to balls hitting each other,  while even collisions, $n=2k$, to the light ball hitting the wall.
Then the time between the consecutive collisions is
\begin{equation}
\begin{split}
\tau_{2k} &= \frac{x_{2k}-X_{2k}}{V_{2k}-v_{2k}} \\
\tau_{2k-1} &=  - \frac{x_{2k-1}}{v_{2k-1}},
\end{split}
\label{Eq:tau(n):exact}
\end{equation}
where $\tau_{2k} = t_{2k+1}- t_{2k}$ is the time interval passed between ball--ball $2k$ and the subsequent ball--wall $2k+1$ collision,
and $\tau_{2k-1} = t_{2k} - t_{2k-1}$ is the time interval between ball--wall $2k-1$ and ball--ball $2k$ collisions.
The time moment of the $n$-th collision can be calculated as the sum of the preceding time intervals as
\begin{equation}
t_n = \sum\limits_{l=1}^{n-1} \tau_l.
\label{Eq:t(n):exact}
\end{equation}

The solution to the equations of motion can be expressed as the following iterative formulas for a ball--ball collision ($n=2k+1, k=0,1,\dots$)
\begin{equation}
\begin{split}
X_{2k+1} &= X_{2k} + (X_{2k}-x_{2k}) \frac{V_{2k}}{v_{2k}-V_{2k}} \\
x_{2k+1} &= x_{2k} + (X_{2k}-x_{2k}) \frac{v_{2k}}{v_{2k}-V_{2k}} \\
V_{2k+1} &= \frac{M-m}{m+M} V_{2k} + \frac{2m}{m+M} v_{2k} \\
v_{2k+1} &= \frac{2M}{m+M} V_{2k}+ \frac{m-M}{m+M} v_{2k},
\end{split}
\label{Eq:equation of motion:BB}
\end{equation}
and for a ball--wall collision ($n=2k, k=1,2,\dots$)
\begin{equation}
\begin{split}
X_{2k} &=X_{2k-1}-\frac{V_{2k-1}}{v_{2k-1}} x_{2k-1} \\
x_{2k} &=0 \\
V_{2k} &=V_{2k-1} \\
v_{2k} &= -v_{2k-1}. \\
\end{split}
\label{Eq:equation of motion:BW}
\end{equation}

The iterative process stops when one of the following equivalent conditions holds: (i) $X_n>0$; (ii) $V_n$ is not monotone and starts decreasing after the point of return; (iii) after a ball--ball collision $v_n<0$, which physically means that the light ball goes to $-\infty$.

\end{document}